\documentclass[a4paper]{amsart}

\usepackage{amsmath,amssymb,amsfonts,exscale,url,doi,graphicx,xcolor}
\usepackage[ruled,vlined,algosection]{algorithm2e}
\usepackage[numbers,sort]{natbib}

\definecolor{color1}{rgb}{0,0,1}
\hypersetup{pdfborder={1 0 0}}
\renewcommand{\ldots}{\ensuremath{\dotsc}}
\newcommand{\comment}[1]{}

\newcommand{\C}{\mathbb{C}}
\newcommand{\R}{\mathbb{R}}

\newcommand{\X}{\mathcal{X}}
\newcommand{\Y}{\mathcal{Y}}

\newcommand{\K}{\mathcal{K}}

\newcommand{\V}{\mathcal{V}}

\newcommand{\Z}{\mathcal{Z}}

\def\la{\lambda}

\def\span{{\rm span}}
\def\dim{{{\rm dim}}}
\def\max{{{\rm max}}}
\def\min{{{\rm min}}}

\newcommand{\diag}{{\rm{diag}}}

\newtheorem{theorem}{Theorem}[section]
\newtheorem{lemma}[theorem]{Lemma}

\newtheorem{remark}{Remark}[section]
\numberwithin{equation}{section}

\begin{document}

\title[Angle-free bounds for restarted block eigensolvers]
{Angle-free cluster robust Ritz value bounds\\
 for restarted block eigensolvers}

\author[M.~Zhou]{Ming Zhou}
\address{Universit\"at Rostock, Institut f\"ur Mathematik, 
        Ulmenstra{\ss}e 69, 18055 Rostock, Germany}
\email{ming.zhou at uni-rostock (dot) de}

\author[A.~V.~Knyazev]{Andrew V.~Knyazev}
\address{Department of Mathematical and Statistical Sciences,
University of Colorado Denver, Denver, CO 80217-3364, USA}
\email{andrew.knyazev at ucdenver (dot) edu}

\author[K.~Neymeyr]{Klaus Neymeyr}
\address{Universit\"at Rostock, Institut f\"ur Mathematik, 
        Ulmenstra{\ss}e 69, 18055 Rostock, Germany}
\email{klaus.neymeyr at uni-rostock (dot) de}


\keywords{Ritz values, subspace eigensolvers, cluster robustness.
\hfill November 2, 2022}

\begin{abstract}
Convergence rates of block iterations for solving eigenvalue problems typically
measure errors of Ritz values approximating eigenvalues. The errors of the Ritz values
are commonly bounded in terms of principal angles between the initial or iterative subspace
and the invariant subspace associated with the target eigenvalues. Such bounds thus
cannot be applied repeatedly as needed for restarted block eigensolvers, since
the left- and right-hand sides of the bounds use different terms. They must be
combined with additional bounds which could cause an overestimation. Alternative
repeatable bounds that are angle-free and depend only on the errors of the Ritz values
have been pioneered for Hermitian eigenvalue problems in \doi{10.1515/rnam.1987.2.5.371}
but only for a single extreme Ritz value. We extend this result to all Ritz values and achieve
robustness for clustered eigenvalues by utilizing nonconsecutive eigenvalues.
Our new bounds cover the restarted block Lanczos method and its modifications
with shift-and-invert and deflation, and are numerically advantageous.
\end{abstract}

\maketitle
\pagestyle{myheadings}
\thispagestyle{plain}

\section{Introduction}\label{intro}

Fundamental vector iterations for solving matrix eigenvalue problems,
e.g., the power method and the Rayleigh quotient gradient method,
need to be implemented in a block form or combined with deflation
in order to compute several eigenvalues and
the associated eigenvectors. The block implementation is particularly
crucial for multiple or clustered eigenvalues since the convergence speed
of a vector iteration usually depends on consecutive eigenvalues.
For instance, the power method applied to a normal matrix $A\in\C^{n \times n}$ 
with the eigenvalue arrangement $|\la_1|\ge\cdots\ge|\la_n|$
has the convergence factor $|\la_2/\la_1|$ for approximating $\la_1$,
or generally $|\la_{i+1}/\la_i|$ for approximating $\la_i$
where deflation is required for $i>1$.
An acceleration up to $|\la_{i+p}/\la_i|$
can be achieved by the block power method using $p$-dimensional
subspaces instead of vectors as iterates.

Many block eigensolvers can be formulated as accelerated versions
of the block power method applied to a matrix function $f(A)$.
Their convergence factors thus depend on certain eigenvalues of $f(A)$.
In particular, the traditional analysis of the block Lanczos method, e.g.,
by Saad \cite[(3.10)]{s1980} utilizes shifted Chebyshev polynomials
as $f(\cdot)$. However, the resulting bounds could 
cause an overestimation for restarted iterations with low-dimensional subspaces,
and combining them with additional bounds is necessary
for investigating multiple outer steps. In contrast, bounds of the concise form
$\varepsilon_{\mathrm{next}}\le\varphi\,\varepsilon_{\mathrm{current}}$,
where $\varepsilon$ and $\varphi<1$ denote errors and convergence factors,
can be applied repeatedly giving direct bounds
like $\varepsilon^{(\ell)}\le\varphi^{\ell}\,\varepsilon^{(0)}$
showing the cumulative decrease factor $\varphi^{\ell}$ of the initial error
after $\ell$ outer steps of a restarted block iteration.

We discover such repeatable bounds
for Ritz value errors by extending convergence analysis
of an abstract block iteration by Knyazev \cite{k,ks} for Hermitian eigenvalue problems.
Therein the error factor $\varphi$ of the eigenvalue approximation by Ritz values is defined by distance ratios
with respect to relevant eigenvalues, and the central part
of the derivation utilizes monotonicity of the Rayleigh quotient
rather than angle-based inequalities. Applying these bounds
to the restarted block Lanczos method \cite{cw} overcomes limitations
of traditional angle-dependent bounds from \cite{s1980,ks,lz2015} and
enhances results \cite{z2018} concerning cluster robustness.
Further applications are considered for shift-and-invert and deflation.

We finally mention that multiplying a subspace by the matrix function $f(A)$
is interpreted in \cite{zakn2022} as applying a graph-based filter to a block of vectors.
For image denoising, this approach can extend the single-vector filtering proposed
in \cite{km2015} to multiple, e.g., multi-spectral, images.

\subsection{Ritz value bounds for an abstract block iteration}

Our investigation of restarted block eigensolvers focuses on the abstract block iteration
\begin{equation}\label{abi}
 \Y'=f(A)\Y
\end{equation}
with successive iterative subspaces $\Y$ and $\Y'$.
The function $f(\cdot)$ can be specified afterwards for concrete eigensolvers.
The original introduction of \eqref{abi} in \cite{k,ks} begins with
a generalized eigenvalue problem. Therein a pair $(M,L)$ of Hermitian matrices
with positive definite $L$ is considered, and \eqref{abi} can be reformulated
for the operator $L^{-1}M$ that is Hermitian with respect to the inner product induced by $L$.
To simplify notation, we use without loss of generality
the theoretically equivalent setup $A=L^{-1/2}ML^{-1/2}$
with respect to the Euclidean inner product so that $A$ is a Hermitian matrix.
The eigenvalue problems of $A$ and $(M,L)$ can easily be converted into each other:
\[AX=X\Lambda \quad\Leftrightarrow\quad
 MV=LV\Lambda \quad\mbox{for}\quad
 X=L^{1/2}V,\]
where $\Lambda$ is a diagonal matrix containing eigenvalues, and
the columns of $X$ or $V$ are the associated eigenvectors of $A$ or $(M,L)$.

We arrange the eigenvalues of $A$ as $\la_1\ge\cdots\ge\la_n$,
and denote by $p$ the block size of \eqref{abi}, i.e., $p=\dim\Y$.
As a single-step iteration, \eqref{abi} describes an outer step of a restarted block iteration
which aims at the $p$ largest eigenvalues of $A$, i.e.,
the subspace $\Y$ is the current approximation
and the subspace $\Y'$ is an improved approximation
to the invariant subspace $\X$ associated with $\la_1,\ldots,\la_p$. This consideration
can be extended to other target eigenvalues after 
elementary reformulations, as introduced below in Section 4
concerning applications of our analysis and concrete eigensolvers.

Under natural assumptions (cf.~Theorem \ref{thm:afrvb}), the subspace $\Y'$
also has dimension $p$, and provides more accurate
Ritz values as approximate eigenvalues in comparison to those by the subspace $\Y$.

We first review an angle-dependent Ritz value bound based on \cite[(2.20)]{ks}.
By using arbitrary orthonormal basis matrices $X$ and $Y$
of the subspaces $\X$ and $\Y$ mentioned above, the Euclidean angle $\angle(\X,\Y)$
is given by $\arccos\tau$ with the smallest singular value $\tau$ of $X^HY$.
Then the $i$th largest Ritz value $\eta'_i$ in $\Y'$
for an arbitrary index $i\in\{1,\ldots,p\}$ is shown in \cite[(2.20)]{ks} to fulfill 
\begin{equation}\label{adrvb}
 (\la_i-\eta'_i)/(\eta'_i-\la_n)
 \le\varphi_i^2\tan^2\angle(\X,\Y)
\end{equation}
with the convergence factor
\begin{equation}\label{cf}
 \varphi_i=\big(\max_{j\in\{p+1,\ldots,n\}}|f(\la_j)|\big)
  /\big(\min_{j\in\{1,\ldots,i\}}|f(\la_j)|\big).
\end{equation}

The derivation of \eqref{adrvb}
uses the intermediate angle bound
\begin{equation}\label{angleb}
 \tan\angle\big(\X_i,f(A)\Y_i\big)
 \le\varphi_i\tan\angle(\X_i,\Y_i)
\end{equation}
for an invariant subspace $\X_i$ associated with $\la_1,\ldots,\la_i$
and an auxiliary subspace $\Y_i\subseteq\Y$
whose orthogonal projection on $\X$ coincides with $\X_i$.
Bound \eqref{angleb} implies \eqref{adrvb} by using
the perturbation bound
\[(\la_i-\widetilde{\eta}_i)/(\widetilde{\eta}_i-\la_n)
 \le\tan^2\angle\big(\X_i,f(A)\Y_i\big)\]
for the $i$th largest Ritz value $\widetilde{\eta}_i$ in $f(A)\Y_i$
together with the inequalities $\la_i\ge\eta'_i\ge\widetilde{\eta}_i$ and
$\angle(\X_i,\Y_i)=\angle(\Y_i,\X)\le\angle(\X,\Y)$.

Our recent paper \cite{zakn2022} upgrades \eqref{adrvb}
by majorization-based techniques from \cite{ka10}. The resulting bounds
are concerned with partial sums of principal angles and Ritz value errors. Therein
tuplewise convergence factors enables improving
comparable bounds from \cite{lz2015} which use scalar convergence factors.

Bound \eqref{adrvb} and its upgrades cannot be applied repeatedly as needed
for restarted block eigensolvers, since the left- and right-hand sides of the bounds
use different terms. Alternative repeatable bounds that are angle-free
and depend only on the errors of the Ritz values have been pioneered
for Hermitian eigenvalue problems in \cite[(2.22)]{ks}
but only for a single extreme Ritz value,
\begin{equation}\label{afrvb}
 (\la_p-\eta'_p)/(\eta'_p-\la_{p+1})
 \le\varphi_p^2\,(\la_p-\eta_p)/(\eta_p-\la_{p+1})
\end{equation}
using the smallest ($p$th largest) Ritz value $\eta_p$ in $\Y$
rather than the angle $\angle(\X,\Y)$. We note that \eqref{afrvb}
additionally requires the assumptions
$\eta_p>\la_{p+1}$ and $|f(\la_1)| \ge\cdots\ge |f(\la_p)|$.
Correspondingly, the convergence factor $\varphi_p$ turns into
$\big(\max_{j\in\{p+1,\ldots,n\}}$ $|f(\la_j)|\big)/|f(\la_p)|$.

Similar bounds for other Ritz values can be derived by
adapting bound \eqref{afrvb} to the partial iteration
$\widetilde{\Y}'_i=f(A)\widetilde{\Y}_i$ of \eqref{abi}
where the subspace $\widetilde{\Y}_i$ is
spanned by Ritz vectors associated with the $i$ largest Ritz values in $\Y$.
This yields
\begin{equation}\label{afrvbi}
 (\la_i-\eta'_i)/(\eta'_i-\la_{i+1})
 \le\big((\max_{j\in\{i+1,\ldots,n\}}|f(\la_j)|)/|f(\la_i)|\big)^2\,(\la_i-\eta_i)/(\eta_i-\la_{i+1})
\end{equation}
with an alternative convergence factor.
A specified form of \eqref{afrvbi} for the block steepest ascent method
for Hermitian eigenvalue problems appears in \cite[(2.5.8)]{k}.
However, the convergence factor in \eqref{afrvbi}
could be close to $1$ for clustered eigenvalues
and in the case $i=1$ simply turns into a convergence factor
for vector iterations.
Thus \eqref{afrvbi} is unable to reflect cluster robustness of block eigensolvers.

Repeatable angle-free bounds involving exclusively errors in Ritz values, such as
\eqref{afrvb} and \eqref{afrvbi}, are naturally used to analyze solvers that optimize
the Ritz values, e.g.,\ the (block) Lanczos method.
We aim to upgrade these bounds in the present work.
Repeatable angle-based bounds for Ritz vectors in such solvers
are apparently only known for single-vector versions; see
\cite[(2.2)]{ks91} for the steepest ascent method
and \cite[(3.3)]{nz2016} for the restarted Lanczos method.
Deriving sharp bounds for Ritz vectors in block solvers
is one of our future topics.

\subsection{Block Lanczos method}

One can interpret iteration \eqref{abi} as an $A$-depend\-ent filter
on the subspace $\Y$ aiming in our case to move $\Y$ closer to
the $A$-invariant subspace $\X$, as can be seen, e.g., from bound \eqref{angleb}.
If the function $f(\cdot)$ is explicitly given and does not depend on $\Y$,
the corresponding convergence factors can be easily explicitly derived.
For example, assuming that $A$ is positive semi-definite,
the filter $f(A)=A^k$ by the $k$-step power method trivially gives
$\varphi_i=\lambda_{p+1}^k/\lambda_i^k$ in \eqref{cf}; cf.~\cite[Section 2]{r1969}.
Similar results for indefinite $A$ are available in \cite[Section V]{s1969} and \cite[Subsection 5.1]{gwbw}.
A basic proof technique therein is to skip several eigenvalues
by using the orthogonal complement of the associated invariant subspace.

Advanced eigensolvers, such as the block Lanczos method, perform an extra step
of determining the output subspace $\Y'$ via Rayleigh-Ritz optimization, in our case
maximizing the Ritz values on $f(A)\Y$ with respect to arbitrary $f(\cdot)$.
While the combined procedure may still be
technically expressed by formula \eqref{abi}, the underlying optimal $f(\cdot)$
now depends on the input subspace $\Y$, i.e., can be viewed as a nonlinear
implicitly defined filter, and a sharp convergence factor $\varphi_i$ may be difficult to determine.
Thus, instead of specifying bounds like \eqref{adrvb} and \eqref{afrvb} directly
with an unknown optimal $f(\cdot)$, one first uses a surrogate function as $f(\cdot)$,
commonly based on Chebyshev polynomials, with explicitly known convergence factors
$\varphi_i$ in \eqref{adrvb} and \eqref{afrvb}.
Then final bounds follow from the optimality of the Ritz values.  

In the block Lanczos method, with any polynomial
of degree $k$ as $f(\cdot)$, the subspace $\Y'$ from \eqref{abi} is a subset
of the block Krylov subspace $\K=\Y+A\Y+\cdots+A^{k}\Y$.
Then the $i$th largest Ritz value $\psi_i$ in $\K$
is not smaller than $\eta'_i$ according to the Courant-Fischer principles.
A typical choice of the surrogate function
to be used for determination of the convergence factors is
\begin{equation}\label{fch}
 f(\alpha)=T_{k}\left(1+2\,\frac{\alpha-\la_{p+1}}{\la_{p+1}-\la_n}\right)
\end{equation}
with the Chebyshev polynomial $T_{k}$ of the first kind.
Then the convergence factor $\varphi_i$
in \eqref{adrvb} and \eqref{afrvb} becomes
\begin{equation}\label{cf1}
 \sigma_i=\left[T_{k}\left(1+2\,
 \frac{\la_i-\la_{p+1}}{\la_{p+1}-\la_n}\right)\right]^{-1}.
\end{equation}
This leads to the following bounds proposed in \cite[Section 2.6]{k}:
\begin{equation}\label{adrvb1}
 (\la_i-\psi_i)/(\psi_i-\la_n)
 \le\sigma_i^2\tan^2\angle(\X,\Y),
\end{equation}
\begin{equation}\label{afrvb1}
 (\la_p-\psi_p)/(\psi_p-\la_{p+1})
 \le\sigma_p^2\,(\la_p-\eta_p)/(\eta_p-\la_{p+1}).
\end{equation}
The majorization-type bounds from \cite{zakn2022}
can be specified in a similar way.

The convergence factor $\sigma_i$ defined in \eqref{cf1}
decreases asymptotically with $k\to\infty$ as a geometric progression
with the reduction rate
\[\frac{\sqrt{\kappa_i}-1}{\sqrt{\kappa_i}+1}\quad\mbox{for}\quad
 \kappa_i=\frac{\la_i-\la_n}{\la_i-\la_{p+1}}.\]
Therefore \eqref{adrvb1} or \eqref{afrvb1}
can predict a rapid convergence of $\psi_i\to\la_i$
or $\psi_p\to\la_p$ provided that the gap between
$\la_p$ and $\la_{p+1}$ is sufficiently large.

In the particular case $k{\,=\,}2$, the block Lanczos method is reduced to one step
of the block steepest ascent method.
The convergence factor $\sigma_i$ becomes (see \cite[Section 2.5]{k})
\[\frac{\kappa_i-1}{\kappa_i+1}
 \ \quad\mbox{or}\quad \
 \frac{1-\xi_i}{1+\xi_i}\quad\mbox{for}\quad
 \xi_i=\frac{\la_i-\la_{p+1}}{\la_i-\la_n}.\]

Bound \eqref{adrvb1} improves the previously known result
\cite[Theorem 6]{s1980} (with different notations of indices) in two regards:
\begin{itemize}
\item[(i)] The left-hand side is not smaller than
the measure $(\la_i-\psi_i)/(\la_i-\la_n)$ used in \cite{s1980}
(which is more suitable for sine-based bounds
rather than tangent-based bounds).
\item[(ii)] The right-hand side does not contain
critical terms which depend on Ritz values in $\K$
and could deteriorate the bound for clustered eigenvalues.
\end{itemize}

Nevertheless, bound (1.9) is less suitable for restarted eigensolvers
in contrast to angle-free bounds like (1.10) 
pioneered in \cite{k,ks} -- the focus of the present investigation. 

\subsection{Restarted block Lanczos method}

In practice, the block Lanczos meth\-od needs to be restarted
for ensuring numerical stability and reasonable storage requirements.
The corresponding block Krylov subspace is extended up to a fixed index $k$,
and then reduced to a subspace spanned by selected Ritz vectors, e.g.,
\begin{equation}\label{rbl}
 \Y^{(\ell+1)}\quad\xleftarrow{\mathrm{RR}[p]}\quad
 \Y^{(\ell)}+A\Y^{(\ell)}+\cdots+A^{k}\Y^{(\ell)}
\end{equation}
where the Rayleigh-Ritz procedure $\mathrm{RR}[p]$
extracts Ritz vectors associated with the $p$ largest Ritz values.
Bounds \eqref{adrvb1} and \eqref{afrvb1} can easily be adapted to
\eqref{rbl} as single-step Ritz value bounds
\begin{equation}\label{adrvb2}
 (\la_i-\psi_i^{(\ell+1)})/(\psi_i^{(\ell+1)}-\la_n)
 \le\sigma_i^2\tan^2\angle(\X,\Y^{(\ell)})
\end{equation}
and
\begin{equation}\label{afrvb2}
 (\la_p-\psi_p^{(\ell+1)})/(\psi_p^{(\ell+1)}-\la_{p+1})
 \le\sigma_p^2\,(\la_p-\psi_p^{(\ell)})/(\psi_p^{(\ell)}-\la_{p+1}).
\end{equation}

We note that \eqref{adrvb2} cannot directly describe the convergence behavior
of the Ritz value sequence $(\psi_i^{(\ell)})_{\ell\in\mathbb{N}}$.
A possible way out is extending \eqref{adrvb2} in multiple steps. This requires
additional bounds for $\angle(\X,\Y^{(\ell+1)})$ in terms of $\psi_i^{(\ell+1)}$.
One such bound following \cite[pp.~382]{ks} is
\[\sin^2\angle(\X,\Y^{(\ell+1)})\le(\la_1-\psi_p^{(\ell+1)})/(\la_1-\la_{p+1})\]
which however could cause an overestimation, especially
in the case $\psi_p^{(\ell+1)}\le\la_{p+1}$.
Similarly, a more accurate approach based on the intermediate bound \eqref{angleb}
requires a reasonable (and rather challenging) comparison of the angles
$\angle\big(\X_i,f(A)\Y_i^{(\ell)}\big)$ and $\angle\big(\X_i,\Y_i^{(\ell+1)}\big)$.
Ritz value bounds from \cite{s1980,lz2015} and
majorization-type bounds from \cite{zakn2022} also depend on angles
and cannot easily be adapted to multiple steps.
The angle-dependence is thus an obstacle for deriving direct
Ritz value bounds for the restarted block Lanczos method
or other similar methods.

In contrast, angle-free bound \eqref{afrvb2} clearly demonstrates
that $(\psi_p^{(\ell)})_{\ell\in\mathbb{N}}$ converges
with respect to the measure $(\la_p-\psi)/(\psi-\la_{p+1})$.
A repeated application of \eqref{afrvb2} leads to the a priori bound
\[(\la_p-\psi_p^{(\ell)})/(\psi_p^{(\ell)}-\la_{p+1})
 \le\sigma_p^{2\ell}\,(\la_p-\psi_p^{(0)})/(\psi_p^{(0)}-\la_{p+1}).\]
Moreover, \eqref{afrvb2} can be generalized to 
$(\psi_i^{(\ell)})_{\ell\in\mathbb{N}}$ for $i\in\{1,\ldots,p\}$
as in \cite[Section 3]{z2018}, namely,
\begin{equation}\label{afrvb2a}
 (\la_t-\psi_i^{(\ell+1)})/(\psi_i^{(\ell+1)}-\la_{t+1})
 \le\widetilde{\sigma}_t^2\,(\la_t-\psi_i^{(\ell)})/(\psi_i^{(\ell)}-\la_{t+1}).
\end{equation}
Therein $\psi_i^{(\ell)}$ is assumed to be located in an arbitrary eigenvalue interval
$(\la_{t+1},\la_t)$ for $t \ge i$, and the convergence factor is defined by
\[\widetilde{\sigma}_t=\left[T_{k}\left(1+2\,
 \frac{\la_t-\la_{t+1}}{\la_{t+1}-\la_n}\right)\right]^{-1}.\]
Bound \eqref{afrvb2a} for $t=i$ simply follows from \eqref{afrvbi}.
Nevertheless, \eqref{afrvb2a} is rather appropriate for well-separated eigenvalues
on account of the gap between the consecutive eigenvalues $\la_t$ and $\la_{t+1}$.
Generalizing \eqref{afrvb2} with convergence factors like \eqref{cf1} is desirable
for interpreting the practically observed cluster robustness of the restarted block Lanczos method.

\subsection{Aim and outline}

In this paper, we generally consider the abstract block iteration \eqref{abi}
and extend the angle-free Ritz value bound \eqref{afrvb} to
arbitrarily located Ritz values. The error of a Ritz value $\psi$ is measured
by $(\la_t-\psi)/(\psi-\la_{t+1})$ or $(\la_{t-p+i}-\psi)/(\psi-\la_{t+1})$.
The former is related to bounds \eqref{afrvbi} and \eqref{afrvb2a} mentioned above,
whereas the latter is suitable for deriving
convergence factors which are independent of the eigenvalues
$\la_{t-p+i+1},\ldots,\la_t$. In the special case $t=p$, these irrelevant eigenvalues
read $\la_{i+1},\ldots,\la_p$, i.e., those skipped in convergence factors
\eqref{cf} and \eqref{cf1}.

Angle-free Ritz value bounds using such measures are also available
in our investigation of block preconditioned gradient-type eigensolvers;
cf.~\cite{zn2019,zn2022}. These methods can be regarded as perturbed versions
of the block power method or a restarted block Krylov subspace iteration
of degree $2$. Some proof techniques in \cite{zn2019,zn2022} are motivated
by analyzing the case of exact shift-inverse preconditioning
in connection with special forms of \eqref{abi}.
The present paper aims at generalizing this approach to \eqref{abi}
as a basis for upgrading the analysis of more efficient 
block preconditioned eigensolvers including LOBPCG \cite{k01}
and the block Davidson method \cite{mrd}.

Preparing for our main analysis, we formulate in Section 2
two important arguments introducing decisive auxiliary vectors and intermediate bounds.
By combining these arguments in Section 3, we get desirable extensions of \eqref{afrvb}.
Therein a partial iteration of \eqref{abi} is constructed
with certain subspaces orthogonal to the invariant subspace
associated with $\la_{i+1},\ldots,\la_p$ for arbitrary $i\in\{1,\ldots,p\}$.
The approximation of $\lambda_i$ by Ritz values from this partial iteration
is analyzed analogously to \eqref{afrvb} with a novel accuracy measure for block eigensolvers
that allows deriving repeatable cluster-friendly bounds extending \eqref{afrvb} to all Ritz values.
Specifically, under the same assumptions as for \eqref{afrvb}
concerning the final phase of a restarted block eigensolver,
we derive in Theorem \ref{thm:afrvb1} the repeatable bound
\[(\la_i-\widetilde{\eta}'_i)/(\widetilde{\eta}'_i-\la_{p+1})
 \le\frac{\max_{j\in\{p+1,\ldots,n\}}|f(\la_j)|^2}{|f(\la_i)|^2}\
 (\la_i-\widetilde{\eta}_i)/(\widetilde{\eta}_i-\la_{p+1})\]
for the $i$th largest Ritz values $\widetilde{\eta}_i$
and $\widetilde{\eta}'_i$ generated by the partial iteration.
Theorem \ref{thm:afrvb2} utilizes a similar partial iteration
within the orthogonal complement of the invariant subspace
associated with $\la_{t-p+i+1},\ldots,\la_t$, and generalizes the above bound
to an arbitrary outer step of a restarted block eigensolver,
\[(\la_{t-p+i}-\widetilde{\eta}'_i)/(\widetilde{\eta}'_i-\la_{t+1})
 \le\frac{\max_{j\in\{t+1,\ldots,n\}}|f(\la_j)|^2}{|f(\la_{t-p+i})|^2}\
 (\la_{t-p+i}-\widetilde{\eta}_i)/(\widetilde{\eta}_i-\la_{t+1}).\]

These repeatable bounds enable explicit multi-step error estimation for various
restarted block eigensolvers. Therein convergence factors 
with nonconsecutive eigenvalues are preserved and
improve bounds like \eqref{afrvbi} in the case of clustered eigenvalues.

Section 4 introduces applications to the restarted block Lanczos method
and its modifications concerning shift-and-invert and deflation.
Section 5 presents numerical examples for demonstrating our new results.
Some detailed proofs are given in Section 6 (Appendix).

Frequently used symbols and settings are collected in Lemma \ref{lm:afrvb1}.

\section{Preliminaries}

A remarkable argument for deriving the angle-free Ritz value bound \eqref{afrvb}
is concerned with inequalities of the Rayleigh quotient where two vectors
are compared with respect to their orthogonal projections on eigenspaces
\cite[Lemma 2.3.2]{k} by Knyazev. Therewith certain auxiliary vectors
can be constructed in low-dimensional subspaces and
produce an appropriate intermediate bound. In Subsection 2.1, we introduce
this argument within a modified derivation of \eqref{afrvb}
for preparing extensions of \eqref{afrvb}.
A further argument introduced in Subsection 2.2
is related to biorthogonal vectors originally used in the analysis
of the block power method by Rutishauser \cite{r1969}.
We construct similar vectors for skipping irrelevant eigenvalues
mentioned in Subsection 1.4.

\subsection{Derivation of a basic angle-free bound}

As the starting point of our analysis, we modify the derivation of
the known bound \eqref{afrvb} in a style that facilitates extensions.
For the reader's convenience, we recall some basic settings.

\begin{lemma}\label{lm:afrvb1}
Consider the abstract block iteration \eqref{abi}, i.e.,
$\Y'=f(A)\Y$ for a Hermitian matrix $A\in\C^{n \times n}$.
The eigenvalues of $A$ are arranged as $\la_1\ge\cdots\ge\la_n$,
and $x_1,\ldots,x_n$ are associated orthonormal eigenvectors.
Let $Y\in\C^{n \times p}$ be an arbitrary basis matrix of the subspace $\Y$.
The Ritz values of $A$ in $\Y$ are arranged as $\eta_1\ge\cdots\ge\eta_p$.
If $\eta_p>\la_{p+1}$, and $f(\la_j)\neq0$ for each $j\in\{1,\ldots,p\}$,
then $f(A)Y$ is a basis matrix of $\Y'$, and $\dim\Y'=p$.
\end{lemma}
\noindent\textit{Proof.} $\to$ Subsection \ref{pr1}.

Subsequently, we reformulate \cite[Lemma 2.3.2]{k}
for deriving intermediate bounds.

\begin{lemma}\label{lm:rqw}
With the settings from Lemma \ref{lm:afrvb1},
denote by $\rho(\cdot)$ the Rayleigh quotient with respect to $A$.
Then the following statements hold for arbitrary vectors
$u,v\in\C^n{\setminus}\{0\}$ where $u$ satisfies $\la_l\ge\rho(u)\ge\la_{l+1}$.

(a) If $|x_j^Hv|\ge|x_j^Hu|\ \forall\ j{\in}\{1,\ldots,l\}$
 and $|x_j^Hv|\le|x_j^Hu|\ \forall\ j{\in}\{l{+}1,\ldots,n\}$,
 then $\rho(v)\ge\rho(u)$.

(b) If $|x_j^Hv|\le|x_j^Hu|\ \forall\ j{\in}\{1,\ldots,l\}$
 and $|x_j^Hv|\ge|x_j^Hu|\ \forall\ j{\in}\{l{+}1,\ldots,n\}$,
 then $\rho(u)\ge\rho(v)$.
\end{lemma}

We omit the proof of Lemma \ref{lm:rqw} and refer the reader
to \cite[Lemma 3.2]{nz2016}. Although the eigenvalues are assumed
to be simple in \cite{k,nz2016}, the proof therein does not utilize strict inequalities
in $\la_l\ge\rho(u)\ge\la_{l+1}$ and is thus compatible with multiple eigenvalues.

Lemma \ref{lm:rqw} motivates some auxiliary vectors by reweighting
projections of a Ritz vector with respect to eigenvectors.

\begin{lemma}\label{lm:afrvb2}
Following Lemma \ref{lm:afrvb1}, an arbitrary Ritz vector $y'$
associated with the smallest ($p$th largest)
Ritz value $\eta'_p$ in $\Y'$ can be represented
by \,$y'=f(A)Yg$\, with a certain $g\in\C^p{\setminus}\{0\}$.
Assume in addition \,$\la_p>\eta'_p$\,
and \,$|f(\la_1)|\ge\cdots\ge|f(\la_p)|>\nu_p$
for $\nu_p=\max_{j\in\{p+1,\ldots,n\}}|f(\la_j)|$.
Denote by $\rho(\cdot)$ the Rayleigh quotient with respect to $A$. Then
\begin{equation}\label{afrvbr}
 \la_p>\eta'_p=\rho(y')\ge\rho(y^{\circ})\ge\rho(y)\ge\eta_p>\la_{p+1}
\end{equation}
holds for \,$y=Yg$\, and \,$y^{\circ}=f(\la_p)\sum_{j=1}^px_jx_j^Hy
 \,+\, \nu_p\sum_{j=p+1}^nx_jx_j^Hy$. Moreover,
\begin{equation}\label{afrvbs}
 \frac{\rho(y_p)-\rho(y^{\circ})}{\rho(y^{\circ})-\rho(y-y_p)}
 =\frac{\nu_p^2}{|f(\la_p)|^2}\,
 \frac{\rho(y_p)-\rho(y)}{\rho(y)-\rho(y-y_p)}
\end{equation}
holds for \,$y_p=\sum_{j=1}^px_jx_j^Hy$.
\end{lemma}
\noindent\textit{Proof.} $\to$ Subsection \ref{app1}.

Now the derivation of \eqref{afrvb} can be completed by combining intermediate bounds.

\begin{theorem}\label{thm:afrvb}
Consider the abstract block iteration \eqref{abi}. The eigenvalues of $A$
are arranged as $\la_1\ge\cdots\ge\la_n$, and the Ritz values of $A$ in $\Y$
as $\eta_1\ge\cdots\ge\eta_p$ for $p=\dim\Y$. If $\eta_p>\la_{p+1}$, and
\,$|f(\la_1)|\ge\cdots\ge|f(\la_p)|>\max_{j\in\{p+1,\ldots,n\}}|f(\la_j)|$, then
$\dim\Y'=p$, and the smallest ($p$th largest) Ritz value $\eta'_p$ in $\Y'$
fulfills \eqref{afrvb}, i.e.,
\[\frac{\la_p-\eta'_p}{\eta'_p-\la_{p+1}}
 \le\frac{\max_{j\in\{p+1,\ldots,n\}}|f(\la_j)|^2}{|f(\la_p)|^2}\
 \frac{\la_p-\eta_p}{\eta_p-\la_{p+1}}.\]
\end{theorem}
\noindent\textit{Proof.} $\to$ Subsection \ref{pr2}.

\begin{remark}\label{rm:afrvb}
The above derivation of \eqref{afrvb}
is based on \cite[pp.~79--83]{k}, but provides more details
which are crucial for extending \eqref{afrvb}
and deriving more flexible angle-free Ritz value bounds.
In our main analysis, this approach is adapted to invariant subspaces
which enclose appropriate subsets of the iterative subspaces
from \eqref{abi}. The orthogonal complements of such invariant subspaces
are associated with irrelevant eigenvalues mentioned in Subsection 1.4.
\end{remark}

\subsection{Skipping eigenvalues}

The auxiliary vectors introduced in Lemma \ref{lm:afrvb2}
are related to the convergence measure $(\la_p-\psi)/(\psi-\la_{p+1})$
with two consecutive eigenvalues. This fact inspires an analogous approach
within the invariant subspace $\span\{x_{i+1},\ldots,x_p\}^{\perp}$.
Therein $\la_i$ and $\la_{p+1}$ become neighboring eigenvalues
so that intermediate bounds like \eqref{afrvbs} enable bounds
for the $i$th largest Ritz value with respect to $(\la_i-\psi)/(\psi-\la_{p+1})$.
For this purpose, we introduce an appropriate subset
of the subspace $\Y$ from \eqref{abi}.

\begin{lemma}\label{lm:aux}
With the settings from Lemma \ref{lm:afrvb1}, the intersection
$\Y_{[1,i]\cup(p,n]}=\X_{[1,i]\cup(p,n]}\cap\Y$ of the invariant subspace
$\X_{[1,i]\cup(p,n]}=\span\{x_{i+1},\ldots,x_p\}^{\perp}$ and $\Y$ 
has at least dimension $i$. If $\eta_p>\la_{p+1}$, then
$\dim\Y_{[1,i]\cup(p,n]}=i$, and the orthogonal projection
of $\Y_{[1,i]\cup(p,n]}$ on the invariant subspace
$\X=\span\{x_1,\ldots,x_p\}$ coincides with $\span\{x_1,\ldots,x_i\}$.
In addition, if $f(\la_j)\neq0$ for each $j\in\{1,\ldots,i\}$,
then $\Y_{[1,i]\cup(p,n]}'=f(A)\Y_{[1,i]\cup(p,n]}$ also has dimension $i$.
\end{lemma}
\noindent\textit{Proof.} $\to$ Subsection \ref{pr3}.

\begin{remark}\label{rm:aux}
Lemma \ref{lm:aux} suggests a partial iteration
$\Y_{[1,i]\cup(p,n]}'{\,=\,}f(A)\Y_{[1,i]\cup(p,n]}$ of \eqref{abi} where the block size is $i$
and the target invariant subspace is $\span\{x_1,\ldots,x_i\}$.
By adapting the derivation of \eqref{afrvb} to this partial iteration, we get a bound
for the smallest ($i$th largest) Ritz value in $\Y_{[1,i]\cup(p,n]}'$.
Subsequent simple transformations lead to
a cluster robust convergence factor for \eqref{abi}
where the eigenvalues $\la_{i+1},\ldots,\la_p$ are skipped.

The property $\dim\Y_{[1,i]\cup(p,n]}=i$ can also be ensured by
assuming $\angle(\X,\Y)<\pi/2$ as in derivations
of angle-dependent bounds \cite{k,lz2015}
and majorization-type bounds \cite{ka10,zakn2022}.
This assumption is evidently equivalent to the invertibility
of $X^HY$ in the proof of Lemma \ref{lm:aux}.
Moreover, the columns of $Y(X^HY)^{-1}$ correspond to auxiliary vectors
utilized in \cite{r1969} for analyzing the block power method.

In an analogous analysis concerning the location $\eta_p>\la_{t+1}$
for a certain $t \ge p$, we use
$\widetilde{\X}=\span\{x_{t-p+i+1},\ldots,x_t\}^{\perp}$
instead of $\X_{[1,i]\cup(p,n]}$.
Then $\dim(\widetilde{\X}\cap\Y) \ge i$ also holds, but the case
$\dim(\widetilde{\X}\cap\Y') > i$ cannot easily be excluded by assumptions
on Ritz values or angles. Correspondingly, we consider
an arbitrary $i$-dimensional subset of $\widetilde{\X}\cap\Y$ in the further analysis.
\end{remark}

\section{Main results}

The angle-free Ritz value bound \eqref{afrvb}
(with the settings from Theorem \ref{thm:afrvb}) deals with
the $p$th largest Ritz value $\eta'_p$ in $\Y'$ located in the eigenvalue interval
$(\la_{p+1},\la_p)$. Our main results are three extensions of \eqref{afrvb}
concerning the $i$th largest Ritz value $\eta'_i$ for an arbitrary $i \le p$.
The corresponding eigenvalue intervals are
$(\la_{p+1},\la_i)$, $(\la_{t+1},\la_{t-p+i})$ for $t \ge p$,
and $(\la_{t+1},\la_t)$ for $t \ge i$, respectively.

\subsection{Bound depending on nonconsecutive eigenvalues in the final phase}

The first extension of \eqref{afrvb} is directly based on Lemma \ref{lm:aux}.
Therein the assumption $\eta_p>\la_{p+1}$ is related to
the final phase of a restarted block eigensolver.
We adapt the derivation of \eqref{afrvb}
to the partial iteration $\Y_{[1,i]\cup(p,n]}'{\,=\,}f(A)\Y_{[1,i]\cup(p,n]}$
within the invariant subspace $\span\{x_{i+1},\ldots,x_p\}^{\perp}$.

Similarly to Lemma \ref{lm:afrvb2}, we begin with a Ritz vector
in $\Y_{[1,i]\cup(p,n]}'$ and construct auxiliary vectors 
for producing intermediate bounds.

\begin{lemma}\label{lm:afrvb3}
Following Lemma \ref{lm:aux}, an arbitrary Ritz vector $\widetilde{y}'$
associated with the smallest ($i$th largest) Ritz value $\widetilde{\eta}'_i$
in $\Y_{[1,i]\cup(p,n]}'$ can be represented by \,$\widetilde{y}'=f(A)\widetilde{y}$\,
with a certain $\widetilde{y}\in\Y_{[1,i]\cup(p,n]}{\setminus}\{0\}$. Assume in addition
\,$\la_i>\widetilde{\eta}'_i$\, and \,$|f(\la_1)|\ge\cdots\ge|f(\la_i)|>\nu_p$
for $\nu_p=\max_{j\in\{p+1,\ldots,n\}}|f(\la_j)|$.
Denote by $\rho(\cdot)$ the Rayleigh quotient with respect to $A$, and define
\[\textstyle \widetilde{y}_i=\sum_{j=1}^ix_jx_j^H\widetilde{y}, \qquad
 \widetilde{y}^{\circ}=f(\la_i)\widetilde{y}_i\,+\, \nu_p (\widetilde{y}-\widetilde{y}_i).\]
Then it holds that
\begin{equation}\label{afrvbr1}
 \la_i>\widetilde{\eta}'_i=\rho(\widetilde{y}')\ge\rho(\widetilde{y}^{\circ})
 \ge\rho(\widetilde{y})\ge\widetilde{\eta}_i>\la_{p+1}
\end{equation}
where $\widetilde{\eta}_i$ is the smallest ($i$th largest)
Ritz value in $\Y_{[1,i]\cup(p,n]}$. Moreover,
\begin{equation}\label{afrvbs1}
 \frac{\rho(\widetilde{y}_i)-\rho(\widetilde{y}^{\circ})}
 {\rho(\widetilde{y}^{\circ})-\rho(\widetilde{y}-\widetilde{y}_i)}
 =\frac{\nu_p^2}{|f(\la_i)|^2}\,
 \frac{\rho(\widetilde{y}_i)-\rho(\widetilde{y})}
 {\rho(\widetilde{y})-\rho(\widetilde{y}-\widetilde{y}_i)}.
\end{equation}
\end{lemma}
\noindent\textit{Proof.} $\to$ Subsection \ref{app2}.

Lemma \ref{lm:afrvb3} enables a generalization of Theorem \ref{thm:afrvb}
which provides angle-free bounds for the smallest Ritz value in $\Y_{[1,i]\cup(p,n]}'$
and the $i$th largest Ritz value in $\Y'$.

\begin{theorem}\label{thm:afrvb1}
Consider the abstract block iteration \eqref{abi}. The eigenvalues of $A$
are arranged as $\la_1\ge\cdots\ge\la_n$,
and $x_1,\ldots,x_n$ are associated orthonormal eigenvectors.
Define $\X_{[1,i]\cup(p,n]}=\span\{x_{i+1},\ldots,x_p\}^{\perp}$
for $p=\dim\Y$ and $i\in$ $\{1,\ldots,p\}$
(including $\X_{[1,i]\cup(p,n]}=\C^n$ for $i=p$).
Let $\Y_{[1,i]\cup(p,n]}=\X_{[1,i]\cup(p,n]}\cap\Y$,
and $\Y_{[1,i]\cup(p,n]}'$ $=f(A)\Y_{[1,i]\cup(p,n]}$.

Assume that
the smallest ($p$th largest) Ritz value $\eta_p$ in $\Y$ is larger than $\la_{p+1}$,
and \,$|f(\la_1)|\ge\cdots\ge|f(\la_p)|>\max_{j\in\{p+1,\ldots,n\}}|f(\la_j)|$.
Then $\dim\Y'=p$, and $\dim\Y_{[1,i]\cup(p,n]}=\dim\Y_{[1,i]\cup(p,n]}'=i$.
Moreover, it holds that
\begin{equation}\label{eq:afrvb1a}
 \frac{\la_i-\widetilde{\eta}'_i}{\widetilde{\eta}'_i-\la_{p+1}}
 \le\frac{\max_{j\in\{p+1,\ldots,n\}}|f(\la_j)|^2}{|f(\la_i)|^2}\
 \frac{\la_i-\widetilde{\eta}_i}{\widetilde{\eta}_i-\la_{p+1}}
\end{equation}
for the smallest ($i$th largest) Ritz values $\widetilde{\eta}_i$
and $\widetilde{\eta}'_i$ in $\Y_{[1,i]\cup(p,n]}$ and $\Y_{[1,i]\cup(p,n]}'$.
Consequently, the $i$th largest Ritz value $\eta'_i$ in $\Y_{[1,i]\cup(p,n]}'$ fulfills
\begin{equation}\label{eq:afrvb1}
 \frac{\la_i-\eta'_i}{\eta'_i-\la_{p+1}}
 \le\frac{\max_{j\in\{p+1,\ldots,n\}}|f(\la_j)|^2}{|f(\la_i)|^2}\
 \frac{\la_i-\eta_p}{\eta_p-\la_{p+1}}.
\end{equation}
\end{theorem}
\begin{proof}
The assumptions $\eta_p>\la_{p+1}$ and
\,$|f(\la_1)|\ge\cdots\ge|f(\la_p)|>\max_{j\in\{p+1,\ldots,n\}}$ $|f(\la_j)|$
ensure the dimension statements
according to Theorem \ref{thm:afrvb} and Lemma \ref{lm:aux}.

Bound \eqref{eq:afrvb1a} is trivial for $\la_i=\widetilde{\eta}'_i$.
In the nontrivial case $\la_i>\widetilde{\eta}'_i$,
Lemma \ref{lm:afrvb3} is applicable. Then \eqref{afrvbr1} leads to
\[\rho(\widetilde{w})\ge\la_i>\widetilde{\eta}'_i=\rho(\widetilde{y}')
 \ge\rho(\widetilde{y}^{\circ})\ge\rho(\widetilde{y})
 \ge\widetilde{\eta}_i>\la_{p+1}\ge\rho(\widetilde{z})\]
for $\widetilde{w}=\widetilde{y}_i$, $\widetilde{z}=\widetilde{y}-\widetilde{y}_i$.
Combining this with \eqref{afrvbs1} implies \eqref{eq:afrvb1a}
analogously to \eqref{eq:afrvb}.

Finally, \eqref{eq:afrvb1a} is extended as \eqref{eq:afrvb1} by using
$\la_i\ge\eta'_i\ge\widetilde{\eta}'_i\ge\widetilde{\eta}_i\ge\eta_p>\la_{p+1}$.
\end{proof}

The abstract block iteration \eqref{abi} and the partial iteration
$\Y_{[1,i]\cup(p,n]}'{\,=\,}f(A)$ $\Y_{[1,i]\cup(p,n]}$ both produce
subsets of the trial subspace of a restarted block eigensolver
so that Theorem \ref{thm:afrvb1} gives upper bounds
for the convergence rates of Ritz values therein.
In particular, Theorem \ref{thm:afrvb1} indicates that
the convergence rate of the $i$th largest Ritz value in the final phase 
does not depend on the possibly clustered eigenvalues
$\la_{i+1},\ldots,\la_p$. 

Bounds \eqref{eq:afrvb1a} and \eqref{eq:afrvb1} coincide
with the known bound \eqref{afrvb} in the special case $i=p$,
and are thus extensions of \eqref{afrvb} to $i \le p$.
They supplement \eqref{afrvb} especially for a block size $p$
with $\la_p\approx\la_{p+1}$ but $\la_i\gg\la_{p+1}$. Therein
the convergence factor $\big(\max_{j\in\{p+1,\ldots,n\}}|f(\la_j)|^2\big)/|f(\la_i)|^2$
can be bounded away from $1$. Correspondingly, it is meaningful
to observe individual residual vectors instead of an entire residual matrix
in the stopping criterion of a restarted block eigensolver.

\subsection{Bound depending on nonconsecutive eigenvalues in general}

In the second extension of \eqref{afrvb}, the assumption $\eta_p>\la_{p+1}$
is to be relaxed as $\eta_p>\la_{t+1}$ for $t \ge p$ so that the resulting bound
is applicable to each outer step of a restarted block eigensolver.
We first generalize Lemma \ref{lm:aux} and Lemma \ref{lm:afrvb3}
where the skipped eigenvalues $\la_{i+1},\ldots,\la_p$ become
$\la_{t-p+i+1},\ldots,\la_t$, i.e., the indices are shifted by $t{\,-\,}p$.

\begin{lemma}\label{lm:afrvb4}
With the settings from Lemma \ref{lm:afrvb1},
denote for $t\in\{p,\ldots,n{\,-\,}1\}$ and $i\in\{1,\ldots,p\}$ the invariant subspace
$\span\{x_{t-p+i+1},\ldots,x_t\}^{\perp}$ by $\widetilde{\X}$.
Then the intersection $\widetilde{\Y}=\widetilde{\X}\cap\Y$ 
has at least dimension $i$. Moreover, if $\eta_p>\la_{t+1}$,
and $f(\la_j)\neq0$ for each $j\in\{1,\ldots,t{\,-\,}p{\,+\,}i\}$, then
$\widetilde{\Y}'_i=f(A)\widetilde{\Y}_i$ has dimension $i$ for an arbitrary
$i$-dimensional subspace $\widetilde{\Y}_i\subseteq\widetilde{\Y}$.

Consequently, an arbitrary Ritz vector $\widetilde{y}'$
associated with the smallest ($i$th largest) Ritz value $\widetilde{\eta}'_i$
in $\widetilde{\Y}'_i$ can be represented by \,$\widetilde{y}'=f(A)\widetilde{y}$\,
with a certain $\widetilde{y}\in\widetilde{\Y}_i{\setminus}\{0\}$. Assume in addition
\,$\la_{t-p+i}>\widetilde{\eta}'_i$\, and \,$|f(\la_1)|\ge\cdots\ge|f(\la_{t-p+i})|>\nu_t$
for $\nu_t=\max_{j\in\{t+1,\ldots,n\}}|f(\la_j)|$.
Denote by $\rho(\cdot)$ the Rayleigh quotient with respect to $A$, and define
\[\textstyle \widetilde{w}=\sum_{j=1}^{t-p+i}x_jx_j^H\widetilde{y}, \qquad
 \widetilde{z}=\sum_{j=t+1}^nx_jx_j^H\widetilde{y},\qquad
 \widetilde{y}^{\circ}=f(\la_{t-p+i})\widetilde{w}\,+\, \nu_t \widetilde{z}.\]
Then it holds that
\begin{equation}\label{afrvbr2}
 \la_{t-p+i}>\widetilde{\eta}'_i=\rho(\widetilde{y}')\ge\rho(\widetilde{y}^{\circ})
 \ge\rho(\widetilde{y})\ge\widetilde{\eta}_i>\la_{t+1}
\end{equation}
where $\widetilde{\eta}_i$ is the smallest ($i$th largest)
Ritz value in $\widetilde{\Y}_i$. Moreover,
\begin{equation}\label{afrvbs2}
 \frac{\rho(\widetilde{w})-\rho(\widetilde{y}^{\circ})}
 {\rho(\widetilde{y}^{\circ})-\rho(\widetilde{z})}
 =\frac{\nu_t^2}{|f(\la_{t-p+i})|^2}\,
 \frac{\rho(\widetilde{w})-\rho(\widetilde{y})}
 {\rho(\widetilde{y})-\rho(\widetilde{z})}.
\end{equation}
\end{lemma}
\noindent\textit{Proof.} $\to$ Subsection \ref{app3}.

The subspaces $\widetilde{\Y}_i$ and $\widetilde{\Y}'_i$
from Lemma \ref{lm:afrvb4} motivate a partial iteration of \eqref{abi}
within $\widetilde{\X}$ due to
$\widetilde{\Y}_i\subseteq\widetilde{\X}\cap\Y$
and $\widetilde{\Y}'_i\subseteq\widetilde{\X}\cap\Y'$.
A corresponding upgrade of Theorem \ref{thm:afrvb1}
provides angle-free bounds for the smallest Ritz value
in $\widetilde{\Y}'_i$ and the $i$th largest Ritz value in $\Y'$
concerning the eigenvalue interval $(\la_{t+1},\la_{t-p+i})$.

\begin{theorem}\label{thm:afrvb2}
Consider the abstract block iteration \eqref{abi}. The eigenvalues of $A$
are arranged as $\la_1\ge\cdots\ge\la_n$,
and $x_1,\ldots,x_n$ are associated orthonormal eigenvectors.
Define $\widetilde{\Y}=\widetilde{\X}\cap\Y$ with the invariant subspace
$\widetilde{\X}=\span\{x_{t-p+i+1},\ldots,x_t\}^{\perp}$
for $p=\dim\Y$, $t\in\{p,\ldots,n{\,-\,}1\}$ and $i\in\{1,\ldots,p\}$
(including $\widetilde{\X}=\C^n$ for $i=p$).

Assume that
the smallest ($p$th largest) Ritz value $\eta_p$ in $\Y$ is larger than $\la_{t+1}$,
and \,$|f(\la_1)|\ge\cdots\ge|f(\la_t)|>\max_{j\in\{t+1,\ldots,n\}}|f(\la_j)|$.
Then $\dim\widetilde{\Y} \ge i$, and $\dim\widetilde{\Y}'_i=i$ holds
for $\widetilde{\Y}'_i=f(A)\widetilde{\Y}_i$ with an arbitrary
$i$-dimensional subspace $\widetilde{\Y}_i\subseteq\widetilde{\Y}$.
Moreover, denote by $\widetilde{\eta}_i$ and $\widetilde{\eta}'_i$
the smallest ($i$th largest) Ritz values in $\widetilde{\Y}_i$ and $\widetilde{\Y}'_i$.
Then either $\widetilde{\eta}'_i\ge\la_{t-p+i}$, or
\begin{equation}\label{eq:afrvb2a}
 0<\frac{\la_{t-p+i}-\widetilde{\eta}'_i}{\widetilde{\eta}'_i-\la_{t+1}}
 \le\frac{\max_{j\in\{t+1,\ldots,n\}}|f(\la_j)|^2}{|f(\la_{t-p+i})|^2}\
 \frac{\la_{t-p+i}-\widetilde{\eta}_i}{\widetilde{\eta}_i-\la_{t+1}}.
\end{equation}
Consequently, $\dim\Y'=p$, and the $i$th largest Ritz value $\eta'_i$ in $\Y'$ fulfills
either $\eta'_i\ge\la_{t-p+i}$ or
\begin{equation}\label{eq:afrvb2}
 0<\frac{\la_{t-p+i}-\eta'_i}{\eta'_i-\la_{t+1}}
 \le\frac{\max_{j\in\{t+1,\ldots,n\}}|f(\la_j)|^2}{|f(\la_{t-p+i})|^2}\
 \frac{\la_{t-p+i}-\eta_p}{\eta_p-\la_{t+1}}.
\end{equation}
\end{theorem}
\begin{proof}
The dimension statements $\dim\widetilde{\Y} \ge i$ and $\dim\widetilde{\Y}'_i=i$
directly follow from Lemma \ref{lm:afrvb4}.

Bound \eqref{eq:afrvb2a} is stated in the nontrivial case
$\la_{t-p+i}>\widetilde{\eta}'_i$ which is also considered in Lemma \ref{lm:afrvb4}.
The intermediate results \eqref{afrvbr2} and \eqref{afrvbs2}
imply \eqref{eq:afrvb2a} analogously to \eqref{eq:afrvb}.

Furthermore, $\dim\Y'=p$ is indeed the special form
of the statement $\dim\widetilde{\Y}'_i=i$ for $i=p$.
Bound \eqref{eq:afrvb2} follows from \eqref{eq:afrvb2a} according to
$\la_{t-p+i}>\eta'_i\ge\widetilde{\eta}'_i\ge\widetilde{\eta}_i\ge\eta_p>\la_{t+1}$.
\end{proof}

Theorem \ref{thm:afrvb2} serves to discuss the global convergence behavior
of a restarted block eigensolver. Provided that the block size $p$
exceeds the size of each eigenvalue cluster, i.e., $\la_{t-p+i}$
is not close to $\la_{t+1}$, the convergence factor
$\big(\max_{j\in\{t+1,\ldots,n\}}|f(\la_j)|^2\big)$ $/|f(\la_{t-p+i})|^2$
with a properly defined $f(\cdot)$, e.g., a shifted Chebyshev polynomial,
can well reflect cluster robustness.

\subsection{Bound depending on consecutive eigenvalues}

The above extensions of \eqref{afrvb} give bounds
\eqref{eq:afrvb1} and \eqref{eq:afrvb2} for the $i$th largest Ritz value $\eta'_i$
in $\Y'$ in terms of the $p$th largest Ritz value $\eta_p$ in $\Y$.
The next extension deals with the relation between $\eta'_i$ and $\eta_i$
with respect to the eigenvalue interval $(\la_{t+1},\la_t)$ for $t \ge i$.
The resulting bound includes \eqref{afrvb} by setting $t=i=p$.
We begin with the following auxiliary terms.

\begin{lemma}\label{lm:afrvb5}
With the settings from Lemma \ref{lm:afrvb1},
let $y_1,\ldots,y_p$ be orthonormal Ritz vectors associated with
the Ritz values $\eta_1\ge\cdots\ge\eta_p$ in $\Y$.
Define $\widetilde{\Y}_i=\span\{y_1,\ldots,y_i\}$ for $i\in\{1,\ldots,p\}$.
If $\eta_i>\la_{t+1}$ for a certain $t\in\{i,\ldots,n{\,-\,}p{\,+\,}i{\,-\,}1\}$,
and $f(\la_j)\neq0$ for each $j\in\{1,\ldots,t\}$, then
$\widetilde{\Y}'_i=f(A)\widetilde{\Y}_i$ has dimension $i$.

Consequently, an arbitrary Ritz vector $\widetilde{y}'$
associated with the smallest ($i$th largest) Ritz value $\widetilde{\eta}'_i$
in $\widetilde{\Y}'_i$ can be represented by \,$\widetilde{y}'=f(A)\widetilde{y}$\,
with a certain $\widetilde{y}\in\widetilde{\Y}_i{\setminus}\{0\}$. Assume in addition
\,$\la_t>\widetilde{\eta}'_i$\, and \,$|f(\la_1)|\ge\cdots\ge|f(\la_t)|>\nu_t$
for $\nu_t=\max_{j\in\{t+1,\ldots,n\}}$ $|f(\la_j)|$.
Denote by $\rho(\cdot)$ the Rayleigh quotient with respect to $A$, and define
\[\textstyle \widetilde{w}=\sum_{j=1}^tx_jx_j^H\widetilde{y}, \qquad
 \widetilde{z}=\sum_{j=t+1}^nx_jx_j^H\widetilde{y},\qquad
 \widetilde{y}^{\circ}=f(\la_t)\widetilde{w}\,+\, \nu_t \widetilde{z}.\]
Then it holds that
\begin{equation}\label{afrvbr3}
 \la_t>\widetilde{\eta}'_i=\rho(\widetilde{y}')\ge\rho(\widetilde{y}^{\circ})
 \ge\rho(\widetilde{y})\ge\eta_i>\la_{t+1},
\end{equation}
\begin{equation}\label{afrvbs3}
 \frac{\rho(\widetilde{w})-\rho(\widetilde{y}^{\circ})}
 {\rho(\widetilde{y}^{\circ})-\rho(\widetilde{z})}
 =\frac{\nu_t^2}{|f(\la_t)|^2}\,
 \frac{\rho(\widetilde{w})-\rho(\widetilde{y})}
 {\rho(\widetilde{y})-\rho(\widetilde{z})}.
\end{equation}
\end{lemma}
\begin{proof}
Despite the different definition of $\widetilde{\Y}_i$,
we can formally reuse the proof of Lemma \ref{lm:afrvb4}
from Subsection \ref{app3}. A slight modification with the substitution
$t{\,-\,}p{\,+\,}i\to t$ proves Lemma \ref{lm:afrvb5}.
\end{proof}

In contrast to Lemma \ref{lm:afrvb4},
the auxiliary subspaces $\widetilde{\Y}_i$ and $\widetilde{\Y}'_i$
from Lemma \ref{lm:afrvb5} do not build a meaningful partial iteration
since $\widetilde{\Y}_i$ is spanned by Ritz vectors in $\Y$,
but $\widetilde{\Y}'_i$ is not necessarily spanned by Ritz vectors in $\Y'$.
For this reason, we omit in the following theorem 
a Ritz value bound like \eqref{eq:afrvb2a} concerning auxiliary subspaces.

\begin{theorem}\label{thm:afrvb3}
Consider the abstract block iteration \eqref{abi}. The eigenvalues of $A$
are arranged as $\la_1\ge\cdots\ge\la_n$. Assume that
the $i$th largest Ritz value $\eta_i$ in $\Y$ is larger
than $\la_{t+1}$ for a certain $t\in\{i,\ldots,n{\,-\,}p{\,+\,}i{\,-\,}1\}$,
and \,$|f(\la_1)|\ge\cdots\ge|f(\la_t)|>\max_{j\in\{t+1,\ldots,n\}}|f(\la_j)|$.
Then $\dim\Y' \ge i$, and the $i$th largest Ritz value $\eta'_i$ in $\Y'$ fulfills
either $\eta'_i\ge\la_t$ or
\begin{equation}\label{eq:afrvb3}
 0<\frac{\la_t-\eta'_i}{\eta'_i-\la_{t+1}}
 \le\frac{\max_{j\in\{t+1,\ldots,n\}}|f(\la_j)|^2}{|f(\la_t)|^2}\
 \frac{\la_t-\eta_i}{\eta_i-\la_{t+1}}.
\end{equation}
\end{theorem}
\begin{proof}
The subspace $\widetilde{\Y}'_i$ defined in Lemma \ref{lm:afrvb5}
is a subset of $\Y'$. Thus $\dim\Y' \ge \dim\widetilde{\Y}'_i = i$.

Bound \eqref{eq:afrvb3} is stated in the nontrivial case $\la_t>\eta'_i$.
Then $\la_t>\eta'_i\ge\widetilde{\eta}'_i$ allows applying
\eqref{afrvbr3} and \eqref{afrvbs3} so that a bound for $\widetilde{\eta}'_i$
is obtained analogously to \eqref{eq:afrvb}.
This implies \eqref{eq:afrvb3} according to
$\la_t>\eta'_i\ge\widetilde{\eta}'_i>\la_{t+1}$.
\end{proof}

The convergence factor in \eqref{eq:afrvb3} uses consecutive eigenvalues
and is thus less appropriate for describing cluster robustness in comparison to
the convergence factors in \eqref{eq:afrvb1} and \eqref{eq:afrvb2}.
However, the term $(\la_t-\eta_i)/(\eta_i-\la_{t+1})$
can lead to a better bound in the first phase of a restarted block eigensolver,
especially if the relevant eigenvalues are not tightly clustered.

\section{Applications to restarted block eigensolvers}

The angle-free Ritz value bounds from Section 3 can typically be applied to
the restarted block Lanczos method by utilizing shifted Chebyshev polynomials
as $f(\cdot)$ in the abstract block iteration \eqref{abi}. Results by this approach
are presented in Subsection 4.1, beginning with bounds for one outer step
which are comparable with some known angle-dependent bounds
from \cite{s1980,ks}. A decisive advantage of our angle-free bounds
is that applying them to multiple outer steps does not require additional bounds
for connecting successive steps, and thus avoids certain overestimations. 
Subsection 4.2 deals with the convergence analysis
of restarted block eigensolvers with shift-and-invert, e.g.,
for computing eigenvalues of a self-adjoint elliptic partial differential operator.
Therein our angle-free bounds can easily be adapted to
a corresponding generalized eigenvalue problem.
A related discussion on deflation is given in Subsection 4.3.

\subsection{Application to the restarted block Lanczos method}

Following Subsection 1.3, we represent the restarted block Lanczos method
by \eqref{rbl}, and observe the Ritz value sequence
$(\psi_i^{(\ell)})_{\ell\in\mathbb{N}}$ where $\psi_i^{(\ell)}$
is the $i$th largest Ritz value in the $\ell$th iterative subspace $\Y^{(\ell)}$.

By regarding $\Y^{(\ell)}$ as $\Y$ in the abstract block iteration \eqref{abi},
one can select a real polynomial of degree $k$ as $f(\cdot)$
so that $\Y'$ is a subset of $\Y^{(\ell+1)}$, and the $i$th largest Ritz value
$\eta'_i$ in $\Y'$ is a lower bound for $\psi_i^{(\ell+1)}$.
Constructing a ``sharp'' $f(\cdot)$ that enables $\eta'_i=\psi_i^{(\ell+1)}$
involves certain interpolating polynomials which cannot easily be
represented in an explicit form; cf.~\cite[Subsection 3.3]{zn2017}.
In contrast to this, the standard approach from \cite{s1980,p,k}
for investigating the Lanczos method weakly minimizes the convergence factor
concerning an eigenvalue interval. For instance, minimizing
$\varphi_p=\big(\max_{j\in\{p+1,\ldots,n\}}|f(\la_j)|\big)/|f(\la_p)|$
in \eqref{afrvb} can be weakened as minimizing
$\big(\max_{\la\in[\la_n,\la_{p+1}]}|f(\la)|\big)/|f(\la_p)|$.
This results in the shifted Chebyshev polynomial \eqref{fch};
cf.~\cite[Lemma 4.1]{zakn2022}. We reuse \eqref{fch}
together with its alternative
\begin{equation}\label{fch1}
 f(\alpha)=T_{k}\left(1+2\,\frac{\alpha-\la_{t+1}}{\la_{t+1}-\la_n}\right)
\end{equation}
with $t \ge p$ or $t \ge i$ for specifying the new bounds from Section 3.

\begin{theorem}\label{thm:rbl}
Consider the restarted block Lanczos method \eqref{rbl} with the block size $p$.
The eigenvalues of $A$ are arranged as $\la_1\ge\cdots\ge\la_n$, and
the Ritz values of $A$ in $\Y^{(\ell)}$ as $\psi_1^{(\ell)}\ge\cdots\ge\psi_p^{(\ell)}$.
Let $T_{k}$ be the Chebyshev polynomial (of the first kind) of degree $k$,
and $i\in\{1,\ldots,p\}$.

(a) If $\psi_p^{(\ell)}>\la_{t+1}$ for a certain $t\in\{p,\ldots,n{\,-\,}1\}$, then
either $\psi_i^{(\ell+1)}\ge\la_{t-p+i}$, or
\begin{equation}\label{eq:rbl1}
 0<\frac{\la_{t-p+i}-\psi_i^{(\ell+1)}}{\psi_i^{(\ell+1)}-\la_{t+1}}
 \le\left[T_{k}\left(1+2\,\frac{\la_{t-p+i}-\la_{t+1}}{\la_{t+1}-\la_n}\right)\right]^{-2}\
 \frac{\la_{t-p+i}-\psi_p^{(\ell)}}{\psi_p^{(\ell)}-\la_{t+1}}.
\end{equation}
In the latter case, consider orthonormal eigenvectors $x_1,\ldots,x_n$
associated with $\la_1,$ $\ldots,\la_n$. Then the subspace
$\widetilde{\Y}=\span\{x_{t-p+i+1},\ldots,x_t\}^{\perp}\cap\Y^{(\ell)}$
has at least dimension $i$. By using the smallest ($i$th largest) Ritz value 
$\widetilde{\eta}_i$ in an arbitrary $i$-dimensional subspace
$\widetilde{\Y}_i\subseteq\widetilde{\Y}$, it holds that
\begin{equation}\label{eq:rbl2}
 0<\frac{\la_{t-p+i}-\psi_i^{(\ell+1)}}{\psi_i^{(\ell+1)}-\la_{t+1}}
 \le\left[T_{k}\left(1+2\,\frac{\la_{t-p+i}-\la_{t+1}}{\la_{t+1}-\la_n}\right)\right]^{-2}\
 \frac{\la_{t-p+i}-\widetilde{\eta}_i}{\widetilde{\eta}_i-\la_{t+1}}.
\end{equation}

(b) If $\psi_i^{(\ell)}>\la_{t+1}$ for a certain $t\in\{i,\ldots,n{\,-\,}p{\,+\,}i{\,-\,}1\}$, then
either $\psi_i^{(\ell+1)}\ge\la_t$, or
\begin{equation}\label{eq:rbl3}
 0<\frac{\la_t-\psi_i^{(\ell+1)}}{\psi_i^{(\ell+1)}-\la_{t+1}}
 \le\left[T_{k}\left(1+2\,\frac{\la_t-\la_{t+1}}{\la_{t+1}-\la_n}\right)\right]^{-2}\
 \frac{\la_t-\psi_i^{(\ell)}}{\psi_i^{(\ell)}-\la_{t+1}}.
\end{equation}
\end{theorem}
\begin{proof}
The statement (a) follows from Theorem \ref{thm:afrvb2}.
We set $\Y=\Y^{(\ell)}$ and define $f(\cdot)$ by \eqref{fch1}.
Then $\Y'\subseteq\Y^{(\ell+1)}$ so that $\eta'_i\le\psi_i^{(\ell+1)}$.
Moreover, the assumptions in Theorem \ref{thm:afrvb2} are fulfilled:
$\eta_p=\psi_p^{(\ell)}>\la_{t+1}$, and
 \,$|f(\la_1)|\ge\cdots\ge|f(\la_t)|>1=\max_{j\in\{t+1,\ldots,n\}}|f(\la_j)|$.
Specifying \eqref{eq:afrvb2} implies \eqref{eq:rbl1} according to
$\la_{t-p+i}>\psi_i^{(\ell+1)}\ge\eta'_i\ge\eta_p=\psi_p^{(\ell)}>\la_{t+1}$.
An analogous specification of \eqref{eq:afrvb2a} leads to \eqref{eq:rbl2}.

Similarly, the statement (b) is proved by Theorem \ref{thm:afrvb3}
where bound \eqref{eq:afrvb3} is specified.
\end{proof}

Theorem \ref{thm:rbl} extends the angle-free Ritz value bound \eqref{afrvb2}
which focuses on the smallest Ritz value and specifies \cite[(2.22)]{ks}.
The statements deal with all Ritz values. In particular, \eqref{eq:rbl1} is comparable with
existing angle-dependent bounds from \cite[Theorem 6]{s1980} and \cite[(2.20)]{ks}.
The angle-independence enables more accurate predictions in the case
of large angle-dependent factors. Furthermore, the limitation of angle-dependent bounds
for investigating multiple outer steps mentioned in Subsection 1.3
can be overcome by a direct generalization of Theorem \ref{thm:rbl}.

\begin{theorem}\label{thm:rblm}
With the settings from Theorem \ref{thm:rbl}, the following statements
hold for $m$ outer steps of the restarted block Lanczos method \eqref{rbl}.

(a) If $\psi_p^{(\ell)}>\la_{t+1}$ for a certain $t\in\{p,\ldots,n{\,-\,}1\}$, then
either $\psi_i^{(\ell+m)}\ge\la_{t-p+i}$, or
\begin{equation}\label{eq:rblm1}
 0<\frac{\la_{t-p+i}-\psi_i^{(\ell+m)}}{\psi_i^{(\ell+m)}-\la_{t+1}}
 \le\left[T_{k}\left(1+2\,\frac{\la_{t-p+i}-\la_{t+1}}{\la_{t+1}-\la_n}\right)\right]^{-2m}\
 \frac{\la_{t-p+i}-\psi_p^{(\ell)}}{\psi_p^{(\ell)}-\la_{t+1}}.
\end{equation}
In the latter case, consider orthonormal eigenvectors $x_1,\ldots,x_n$
associated with $\la_1,$ $\ldots,\la_n$. Then the subspace
$\widetilde{\Y}=\span\{x_{t-p+i+1},\ldots,x_t\}^{\perp}\cap\Y^{(\ell)}$
has at least dimension $i$. By using the smallest ($i$th largest) Ritz value 
$\widetilde{\eta}_i$ in an arbitrary $i$-dimensional subspace
$\widetilde{\Y}_i\subseteq\widetilde{\Y}$, it holds that
\begin{equation}\label{eq:rblm2}
 0<\frac{\la_{t-p+i}-\psi_i^{(\ell+m)}}{\psi_i^{(\ell+m)}-\la_{t+1}}
 \le\left[T_{k}\left(1+2\,\frac{\la_{t-p+i}-\la_{t+1}}{\la_{t+1}-\la_n}\right)\right]^{-2m}\
 \frac{\la_{t-p+i}-\widetilde{\eta}_i}{\widetilde{\eta}_i-\la_{t+1}}.
\end{equation}

(b) If $\psi_i^{(\ell)}>\la_{t+1}$ for a certain $t\in\{i,\ldots,n{\,-\,}p{\,+\,}i{\,-\,}1\}$, then
either $\psi_i^{(\ell+m)}\ge\la_t$, or
\begin{equation}\label{eq:rblm3}
 0<\frac{\la_t-\psi_i^{(\ell+m)}}{\psi_i^{(\ell+m)}-\la_{t+1}}
 \le\left[T_{k}\left(1+2\,\frac{\la_t-\la_{t+1}}{\la_{t+1}-\la_n}\right)\right]^{-2m}\
 \frac{\la_t-\psi_i^{(\ell)}}{\psi_i^{(\ell)}-\la_{t+1}}.
\end{equation}
\end{theorem}
\begin{proof}
For proving (a), we use the indexed form
\begin{equation}\label{abif}
 \Y^{(\ell+1)}=f(A)\Y^{(\ell)}
\end{equation}
of the abstract block iteration \eqref{abi}, and define $f(\cdot)$ by \eqref{fch1}. 
According to the Courant-Fischer principles, \eqref{rbl}
is a stepwisely accelerated version of \eqref{abif} with respect to Ritz values.
By setting the $\ell$th iterative subspace of \eqref{rbl} as $\Y^{(\ell)}$ in \eqref{abif},
we only need to verify \eqref{eq:rblm1} for \eqref{abif}.
Therein we define $\widetilde{\X}=\span\{x_{t-p+i+1},\ldots,x_t\}^{\perp}$
and $\widetilde{\Y}^{(\ell)}=\widetilde{\X}\cap\Y^{(\ell)}$.
The assumption $\psi_p^{(\ell)}>\la_{t+1}$ and the definition of $f(\cdot)$
allow applying Theorem \ref{thm:afrvb2} to $\Y=\Y^{(\ell)}$
which first shows $\dim\big(f(A)\widetilde{\Y}_i^{(\ell)}\big)=i$
for an arbitrary $i$-dimensional subspace
$\widetilde{\Y}_i^{(\ell)}\subseteq\widetilde{\Y}^{(\ell)}$. Combining this with
\[f(A)\widetilde{\Y}_i^{(\ell)}
 \subseteq\big(f(A)\widetilde{\X}\big)\cap\big(f(A)\Y^{(\ell)}\big)
 \subseteq\widetilde{\X}\cap\Y^{(\ell+1)}=\widetilde{\Y}^{(\ell+1)}\]
ensures that $\widetilde{\Y}_i^{(\ell+1)}=f(A)\widetilde{\Y}_i^{(\ell)}$
is an $i$-dimensional subspace within $\widetilde{\Y}^{(\ell+1)}$
and defines a partial iteration of \eqref{abif}.
In the nontrivial case $\la_{t-p+i}>\psi_i^{(\ell+m)}$ for 
the restarted block Lanczos method \eqref{rbl},
we get $\la_{t-p+i}>\psi_i^{(\ell+m)}\ge\cdots\ge\psi_i^{(\ell+1)}$
so that $\la_{t-p+i}$ is also larger than the corresponding $i$th Ritz values
produced by \eqref{abif} and the partial iteration
$\widetilde{\Y}_i^{(\ell+1)}=f(A)\widetilde{\Y}_i^{(\ell)}$.
Then adapting \eqref{eq:afrvb2a} to the $i$th largest Ritz values
$\widetilde{\eta}_i^{(\ell)}$ and $\widetilde{\eta}_i^{(\ell+1)}$
in $\widetilde{\Y}_i^{(\ell)}$ and $\widetilde{\Y}_i^{(\ell+1)}$ yields
\[0<\frac{\la_{t-p+i}-\widetilde{\eta}_i^{(\ell+1)}}{\widetilde{\eta}_i^{(\ell+1)}-\la_{t+1}}
 \le\left[T_{k}\left(1+2\,\frac{\la_{t-p+i}-\la_{t+1}}{\la_{t+1}-\la_n}\right)\right]^{-2}\
 \frac{\la_{t-p+i}-\widetilde{\eta}_i^{(\ell)}}{\widetilde{\eta}_i^{(\ell)}-\la_{t+1}}.\]
A repeated application thereof leads to a multi-step bound concerning
$\widetilde{\eta}_i^{(\ell)}$ and $\widetilde{\eta}_i^{(\ell+m)}$.
Then \eqref{eq:rblm1} and \eqref{eq:rblm2}
with $\widetilde{\eta}_i=\widetilde{\eta}_i^{(\ell)}$ are verified by considering
$\la_{t-p+i}>\psi_i^{(\ell+m)}\ge\widetilde{\eta}_i^{(\ell+m)}
\ge\widetilde{\eta}_i^{(\ell)}\ge\psi_p^{(\ell)}>\la_{t+1}$.

The statement (b) simply follows from Theorem \ref{thm:rbl} (b)
by repeatedly applying \eqref{eq:rbl3}
in the nontrivial case $\la_t>\psi_i^{(\ell+m)}$.
\end{proof}

Theorem \ref{thm:rblm} indicates that
the Ritz value sequence $(\psi_i^{(\ell)})_{\ell\in\mathbb{N}}$
approaches or exceeds the right end of the considered
eigenvalue interval after sufficiently many outer steps.
Therein the convergence toward an interior eigenvalue $\la_s$
with $s>i$ is not excluded, although rare in practice.
In general, a strict increase $\psi_i^{(\ell+1)}>\psi_i^{(\ell)}$ occurs
if the current iterative subspace $\Y^{(\ell)}$ contains no eigenvectors;
cf.~\cite[Corollary 2]{z2018} or a self-contained explanation
in Subsection \ref{app4}. Based on this fact, the case $\psi_p^{(\ell)}=\la_{t+1}$ 
or $\psi_i^{(\ell)}=\la_{t+1}$ which is not included in the assumption
of Theorem \ref{thm:rblm} can be discussed as follows: If some Ritz vectors
in $\Y^{(\ell)}$ are already eigenvectors, we can adapt Theorem \ref{thm:rblm}
to a reduced iterative subspace spanned by other Ritz vectors.
Otherwise $\psi_p^{(\ell+1)}>\psi_p^{(\ell)}=\la_{t+1}$
or $\psi_i^{(\ell+1)}>\psi_i^{(\ell)}=\la_{t+1}$ holds so that
Theorem \ref{thm:rblm} is applicable after updating the index $t$.

\subsection{Shift-and-invert}

The analysis in Subsection 4.1 is concerned with a Hermitian matrix
$A$ and the computation of its largest eigenvalues. A trivial extension to computing
the smallest eigenvalues can be made by the substitution $A \to -A$.
Now we formulate further extensions starting with
a generalized eigenvalue problem $Lv = \alpha Sv$
for Hermitian matrices $L,S\in\C^{n \times n}$ where $S$ is positive definite.

A basic eigensolver for computing eigenvalues of $(L,S)$ close to
a noneigenvalue shift $\beta\in\R$ is the block shift-and-invert iteration
\begin{equation}\label{bsi}
 \Z^{(\ell+1)}=L_{\beta}^{-1}S\Z^{(\ell)}
 \quad\mbox{with}\quad L_{\beta}=L-\beta S
\end{equation}
which is implemented by solving linear systems for $L_{\beta}$.

In particular, if $\beta$ is smaller than the smallest eigenvalue of $(L,S)$,
then the shifted matrix $L_{\beta}$ is positive definite so that
\eqref{bsi} can be reformulated as
\[L_{\beta}^{1/2}\Z^{(\ell+1)}
 =\big(L_{\beta}^{-1/2}SL_{\beta}^{-1/2}\big)L_{\beta}^{1/2}\Z^{(\ell)}\]
which corresponds to the block power method
for $A=L_{\beta}^{-1/2}SL_{\beta}^{-1/2}$.
Moreover, the eigenvalues of $(L,S)$ and $A$, arranged as
$\alpha_1\le\cdots\le\alpha_n$ and $\la_1\ge\cdots\ge\la_n$, can be
converted into each other by $\la_i=(\alpha_i-\beta)^{-1}$. The conversion between
Ritz values is similar: the $i$th smallest Ritz value $\vartheta_i^{(\ell)}$
of $(L,S)$ in $\Z^{(\ell)}$ and the $i$th largest Ritz value $\psi_i^{(\ell)}$
of $A$ in $\Y^{(\ell)}=L_{\beta}^{1/2}\Z^{(\ell)}$
fulfill $\psi_i^{(\ell)}=(\vartheta_i^{(\ell)}-\beta)^{-1}$.
Therewith bounds from Section 3 can be
specified by $f(\alpha)=\alpha$ and then transformed for \eqref{bsi}.
For instance, \eqref{eq:afrvb1a} leads to a counterpart of \eqref{eq:rblm1}
with the convergence factor $(\la_{p+1}/\la_i)^{2m}$ which is equivalent to
\begin{equation}\label{bsib}
 \frac{\vartheta_i^{(\ell+m)}-\alpha_i}{\alpha_{p+1}-\vartheta_i^{(\ell+m)}}
 \le\left(\frac{\alpha_i-\beta}{\alpha_{p+1}-\beta}\right)^{2m}\
 \frac{\vartheta_p^{(\ell)}-\alpha_i}{\alpha_{p+1}-\vartheta_p^{(\ell)}}.
\end{equation}

Furthermore, a restarted version of the block Davidson method
with $L_{\beta}^{-1}$ as preconditioner can be represented by
\begin{equation}\label{irbl}
 \Z^{(\ell+1)}\quad\xleftarrow{\mathrm{RR}[L,S,p]}\quad
 \Z^{(\ell)}+L_{\beta}^{-1}S\Z^{(\ell)}+\cdots+(L_{\beta}^{-1}S)^{k}\Z^{(\ell)}
\end{equation}
where the Rayleigh-Ritz procedure $\mathrm{RR}[L,S,p]$
extracts Ritz vectors associated with the $p$ smallest Ritz values of $(L,S)$.
As \eqref{irbl} is equivalent to the restarted block Lanczos method \eqref{rbl}
for $A=L_{\beta}^{-1/2}SL_{\beta}^{-1/2}$
and $\Y^{(\ell)}=L_{\beta}^{1/2}\Z^{(\ell)}$,
the specified bounds from Subsection 4.1 can be adapted to \eqref{irbl}
analogously to \eqref{bsib}. For instance, \eqref{eq:rblm1} with $t{\,=\,}p$ corresponds to
\begin{equation}\label{irblb}
 \frac{\vartheta_i^{(\ell+m)}-\alpha_i}{\alpha_{p+1}-\vartheta_i^{(\ell+m)}}
 \le\left[T_{k}\left(1+2\,\frac{(\alpha_i-\beta)^{-1}-(\alpha_{p+1}-\beta)^{-1}}
 {(\alpha_{p+1}-\beta)^{-1}-(\alpha_n-\beta)^{-1}}\right)\right]^{-2m}\
 \frac{\vartheta_p^{(\ell)}-\alpha_i}{\alpha_{p+1}-\vartheta_p^{(\ell)}}.
\end{equation}

The above extension concerning the case $\beta<\alpha_1$
can easily be modified for the case $\beta>\alpha_n$
where $-L_{\beta}$ is positive definite.

In the remaining case $\alpha_1<\beta<\alpha_n$, computing
the largest/smallest eigenvalues smaller/larger than $\beta$
is equivalent to computing the largest eigenvalues of
$(-L_{\beta},M)$ or $(L_{\beta},M)$ for $M=L_{\beta}S^{-1}L_{\beta}$
which is positive definite. This allows again applying bounds from Section 3
where $A$ is defined by $M^{-1/2}(\pm L_{\beta})M^{-1/2}$, i.e.,
a counterpart of $L_{\beta}^{-1/2}SL_{\beta}^{-1/2}$ from the above extension.

\subsection{Deflation}

Our angle-free Ritz value bounds reflect different convergence rates of individual Ritz values.
In practice, the convergence is typically checked by easily computable residual norms
of Ritz vectors. Sufficiently accurate Ritz vectors form a matrix against which
further iterative subspaces are orthogonalized explicitly or implicitly
within a Rayleigh-Ritz procedure. Moreover, one can add random vectors
to the next iterative subspace for keeping the block size unchanged. 
By ignoring errors of accepted Ritz vectors, we can assume that they span
an invariant subspace $\V$, and restrict the investigation
of further steps to the orthogonal complement of $\V$.

As an example, we consider again the generalized eigenvalue problem
from Subsection 4.2, and denote by $v_1,\ldots,v_n$ $S$-orthonormal eigenvectors
associated with the eigenvalues $\alpha_1\le\cdots\le\alpha_n$ of $(L,S)$.
We observe the restarted block eigensolvers \eqref{bsi} and \eqref{irbl}
in the case that the invariant subspace $\V=\span\{v_1,\ldots,v_c\}$ is known.
By deflation, further iterative subspaces are contained in
$\widetilde{\V}=\span\{v_{c+1},\ldots,v_n\}$.

We generally discuss Ritz values in an arbitrary subspace
$\Z\subseteq\widetilde{\V}$ with an $S$-orthonormal basis matrix $Z$ of $\Z$.
The Ritz values of $(L,S)$ in $\Z$ are thus given by the eigenvalues of $Z^HLZ$.
Moreover, by using the matrix $\widetilde{V}=[v_{c+1},\ldots,v_n]$,
we get the $S$-orthogonal projector $P=\widetilde{V}\widetilde{V}^HS$
on $\widetilde{\V}$ so that $Z=PZ$, and
\[Z^HLZ=Z^HP^HLPZ=Z^HS\widetilde{V}\widetilde{V}^HL\widetilde{V}\widetilde{V}^HSZ
 =\widetilde{Z}^HD\widetilde{Z}\]
for $D=\widetilde{V}^HL\widetilde{V}=\diag(\alpha_{c+1},\ldots,\alpha_n)$
and $\widetilde{Z}=\widetilde{V}^HSZ$.
Therein $\widetilde{Z}$ is Euclidean orthonormal since
$\widetilde{Z}^H\widetilde{Z}=Z^HS\widetilde{V}\widetilde{V}^HSZ
 =Z^HSPZ=Z^HSZ$ and $Z$ is $S$-orthonormal.
Thus the eigenvalues of $\widetilde{Z}^HD\widetilde{Z}$
coincide with the Ritz values of $D$ in $\span\{\widetilde{Z}\}$.
In summary, the Ritz values of $(L,S)$ in $\Z$
are just those of $D$ in the subspace $\span\{\widetilde{Z}\}=\widetilde{V}^HS\Z$.

Consequently, an eigensolver for $(L,S)$ after deflation can be represented
by an eigensolver for $D$. We first observe \eqref{bsi}.
The above paragraph shows that the iterate
$\Z^{(\ell)}\subseteq\widetilde{\V}$ shares Ritz values with
$\widetilde{\Z}^{(\ell)}=\widetilde{V}^HS\Z^{(\ell)}$.
A corresponding representation of $\Z^{(\ell)}$ is
$\Z^{(\ell)}=P\Z^{(\ell)}=\widetilde{V}\widetilde{V}^HS\Z^{(\ell)}
 =\widetilde{V}\widetilde{\Z}^{(\ell)}$.
The next iterate $\Z^{(\ell+1)}=L_{\beta}^{-1}S\Z^{(\ell)}$
is also a subset of $\widetilde{\V}$ due to
$L_{\beta}^{-1}S\Z^{(\ell)}\subseteq
 L_{\beta}^{-1}S\widetilde{\V}\subseteq\widetilde{\V}$.
In addition, it holds that
\begin{equation*}\begin{split}
 L_{\beta}^{-1}S\Z^{(\ell)}
 &=PL_{\beta}^{-1}S\Z^{(\ell)}=PL_{\beta}^{-1}SP\Z^{(\ell)}
 =\widetilde{V}\widetilde{V}^HSL_{\beta}^{-1}S
 \widetilde{V}\widetilde{V}^HS\Z^{(\ell)}\\[1ex]
 &=\widetilde{V}\widetilde{V}^HS
 \widetilde{V}D_{\beta}^{-1}\widetilde{V}^HS\Z^{(\ell)}
 =\widetilde{V}D_{\beta}^{-1}\widetilde{\Z}^{(\ell)}
\end{split}\end{equation*}
with $D_{\beta}=\diag(\alpha_{c+1}{\,-\,}\beta,\ldots,\alpha_n{\,-\,}\beta)$.
Thus the subspace $\widetilde{\Z}^{(\ell+1)}=\widetilde{V}^HS\Z^{(\ell+1)}$
representing the next iterate $\Z^{(\ell+1)}$ fulfills
\[\widetilde{\Z}^{(\ell+1)}=\widetilde{V}^HS(L_{\beta}^{-1}S\Z^{(\ell)})
 =\widetilde{V}^HS(\widetilde{V}D_{\beta}^{-1}\widetilde{\Z}^{(\ell)})
 =D_{\beta}^{-1}\widetilde{\Z}^{(\ell)},\]
i.e., \eqref{bsi} is represented by
$\widetilde{\Z}^{(\ell+1)}=D_{\beta}^{-1}\widetilde{\Z}^{(\ell)}$.

Similarly, $(L_{\beta}^{-1}S)^j\Z^{(\ell)}
 =\widetilde{V}D_{\beta}^{-j}\widetilde{\Z}^{(\ell)}$
holds for $j \ge 1$ and gives the representation
\[\widetilde{\Z}^{(\ell+1)}\quad\xleftarrow{\mathrm{RR}[D,p]}\quad
 \widetilde{\Z}^{(\ell)}+D_{\beta}^{-1}\widetilde{\Z}^{(\ell)}
 +\cdots+D_{\beta}^{-k}\widetilde{\Z}^{(\ell)}\]
of \eqref{irbl}. Specifying bounds from Section 3 for these eigensolvers for $D$
produces explicit bounds analogously to Subsection 4.2.
The relevant eigenvalues are from the set
$\{\alpha_{c+1},\ldots,\alpha_n\}$, e.g., a counterpart of \eqref{irblb} reads
\[\frac{\vartheta_i^{(\ell+m)}-\alpha_{c+i}}{\alpha_{c+p+1}-\vartheta_i^{(\ell+m)}}
 \le\left[T_{k}\left(1+2\,
 \mbox{\small$\dfrac{(\alpha_{c+i}-\beta)^{-1}-(\alpha_{c+p+1}-\beta)^{-1}}
 {(\alpha_{c+p+1}-\beta)^{-1}-(\alpha_n-\beta)^{-1}}$}\right)\right]^{-2m}\
 \frac{\vartheta_p^{(\ell)}-\alpha_{c+i}}{\alpha_{c+p+1}-\vartheta_p^{(\ell)}}.\]
Therein the convergence factor can be refined by enlarging the shift $\beta$
up to $\alpha_c$. This enables an acceleration with respect to the number of steps,
but solving linear systems for $L_{\beta}$ with enlarged $\beta$ could be more costly
so that the total computational time is not necessarily reduced.

\section{Numerical examples}

We compare our angle-free Ritz value bounds with their angle-dependent
counterparts which are related to \cite[(2.20)]{ks} by Knyazev and more accurate
than similar traditional bounds from \cite[Theorem 6]{s1980} by Saad
especially for clustered eigenvalues; cf.~the comparison
in \cite[Example 3]{z2018}. The angle-dependent Ritz value bounds
in unitarily invariant norms from \cite[Theorem 8.2]{lz2015} are not included
as they deal with a tuple of Ritz value errors and cannot individually
be applied to the $i$th Ritz value unless $i{\,=\,}1$.

\subsection{Example 1}

We reuse the test matrix from \cite[Example 1]{zakn2022}, i.e.,
the diagonal matrix $A=\mbox{diag}(\la_1,\ldots,\la_n)$ with $n=900$ and
\[\la_1=2, \quad \la_2=1.6, \quad \la_3=1.4, \quad 
 \la_j=1-(j-3)/n \ \ \mbox{for} \ \ j=4,\ldots,n.\]
following \cite[Subsection 4.2]{s1980} and \cite[Example 7.3]{lz2015}.

\begin{figure}[htbp]
\begin{center}
\includegraphics[width=\textwidth]{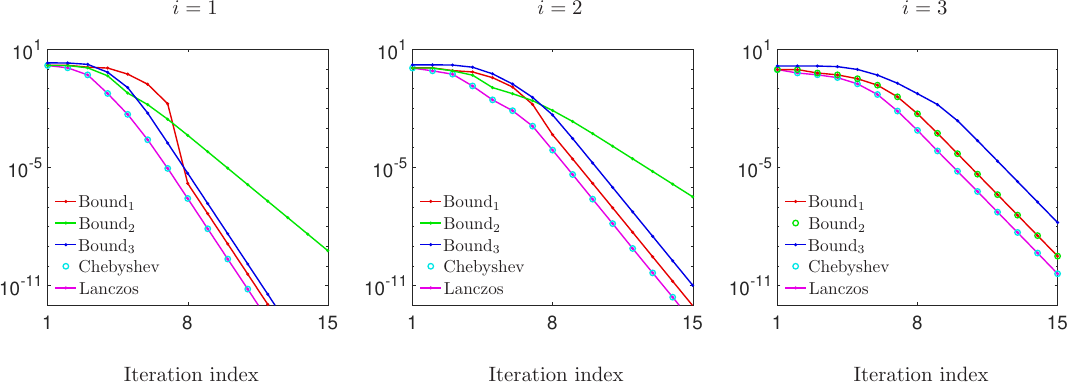}\quad
\end{center}
\par\vskip -2ex
\caption{\small Numerical comparison between several Ritz value bounds
concerning one outer step of the restarted block Lanczos method \eqref{rbl}
in Example 1. The three bound curves are determined by bounds
\eqref{eq:rbl2}, \eqref{eq:rbl3} and \eqref{adrvb3}, respectively.}
\label{fig1a}
\end{figure}

In Figure \ref{fig1a}, we demonstrate bounds from Theorem \ref{thm:rbl}
concerning one outer step of the restarted block Lanczos method \eqref{rbl}.
Therein \eqref{rbl} is implemented with $15$ inner steps and the block size $3$
in each of $1000$ runs with randomly constructed $\Y^{(\ell)}$
and full orthogonalization. We document
the Ritz value error $\la_i-\psi_i^{(\ell+1)}$ for $i\in\{1,\,2,\,3\}$
and each inner step. The associated mean values among $1000$ samples
are displayed by ``Lanczos'' curves. In addition, mean values of Ritz value errors
for $f(A)\Y^{(\ell)}$ with the shifted Chebyshev polynomial \eqref{fch}
are drawn as ``Chebyshev'' circles. We note that the differences
between these two types of data concerning individual Ritz values are less obvious
than comparing Ritz value sums as in \cite[Example 1]{zakn2022}.

Subsequently, we apply bounds \eqref{eq:rbl2} and \eqref{eq:rbl3}
to the first inner steps by determining the index $t$ for the assumption
$\psi_p^{(\ell)}>\la_{t+1}$ or $\psi_i^{(\ell)}>\la_{t+1}$.
Once $\psi_p^{(\ell+1)}$ or $\psi_i^{(\ell+1)}$ in a $(c{\,+\,}1)$th inner step exceeds
an eigenvalue larger than $\la_{t+1}$, we update the index $t$
and observe a subspace $\widehat{\Y}$ spanned by $p$ orthonormal Ritz vectors
associated with the $p$ largest Ritz values in the current block Krylov subspace
$\Y^{(\ell)}+A\Y^{(\ell)}+\cdots+A^{c}\Y^{(\ell)}$. Then the block Krylov subspace
$\widehat{\K}=\widehat{\Y}+A\widehat{\Y}+\cdots+A^{k-c}\widehat{\Y}$
is a subset of $\K=\Y^{(\ell)}+A\Y^{(\ell)}+\cdots+A^{k}\Y^{(\ell)}$.
Thus adapting \eqref{eq:rbl2} and \eqref{eq:rbl3} to $\widehat{\K}$
provides appropriate bounds for $\K$,
i.e., for further inner steps up to the next update.
We convert the evaluated bounds into upper bounds of $\la_i-\psi_i^{(\ell+1)}$.
The associated mean values are displayed in Figure \ref{fig1a}
by ``$\mathrm{Bound}_1$'' and ``$\mathrm{Bound}_2$''
corresponding to \eqref{eq:rbl2} and \eqref{eq:rbl3}.
Furthermore, by using orthonormal eigenvectors $x_1,\ldots,x_p$
associated with $\la_1,\ldots,\la_p$, we evaluate the angle-dependent bound
\begin{equation}\label{adrvb3}
 \frac{\la_i-\psi_i^{(\ell+1)}}{\psi_i^{(\ell+1)}-\la_n}
 \le\left[T_{k}\left(1+2\,\frac{\la_i-\la_{p+1}}{\la_{p+1}-\la_n}\right)\right]^{-2}
 \tan^2\angle(\X_i,\Y_i^{(\ell)})
\end{equation}
in terms of $\X_i=\span\{x_1,\ldots,x_i\}$
and $\Y_i^{(\ell)}=\span\{x_{i+1},\ldots,x_p\}^{\perp}\cap\Y^{(\ell)}$
which can be derived based on \eqref{angleb}
and improves \eqref{adrvb2} (specification of \cite[(2.20)]{ks}) for $i<p$.
Upper bounds of $\la_i-\psi_i^{(\ell+1)}$ generated by \eqref{adrvb3}
are displayed by ``$\mathrm{Bound}_3$''.

The comparison in Figure \ref{fig1a} indicates that $\mathrm{Bound}_1$
is generally more advantageous than the other two bounds.
The overestimation by $\mathrm{Bound}_1$
in several of the first inner steps for $i{\,=\,}1$ is related to the assumption
$\psi_p^{(\ell)}>\la_{t+1}$ of \eqref{eq:rbl2} and the considerably different
convergence rates of $\psi_1^{(\ell)}$ and $\psi_p^{(\ell)}$.
In contrast, $\mathrm{Bound}_2$ using the assumption $\psi_i^{(\ell)}>\la_{t+1}$
provides better alternatives for these steps, but is less accurate
in further steps due to the convergence factor in \eqref{eq:rbl3}
with consecutive eigenvalues. The benefit of $\mathrm{Bound}_3$
is visible in two inner steps for $i{\,=\,}1$. Afterwards the updated
$\mathrm{Bound}_1$ shares the convergence factor with $\mathrm{Bound}_3$,
cf.~\eqref{eq:rbl2} for $t{\,=\,}p$ and \eqref{adrvb3},
and becomes more accurate thanks to the angle-independence.
Moreover, $\mathrm{Bound}_2$ coincides with $\mathrm{Bound}_1$ for $i{\,=\,}3$.

\begin{figure}[htbp]
\begin{center}
\includegraphics[width=\textwidth]{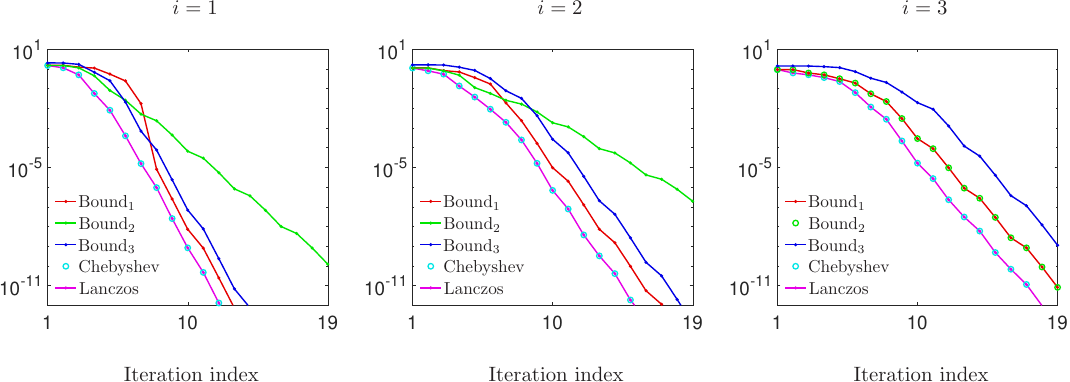}\quad
\end{center}
\par\vskip -2ex
\caption{\small Numerical comparison between several Ritz value bounds
concerning $6$ outer steps ($19$ inner steps)
of the restarted block Lanczos method \eqref{rbl}
in Example 1. The three bound curves are determined by bound combinations
``\eqref{eq:rblm2}$+$\eqref{eq:rbl2}'', ``\eqref{eq:rblm3}$+$\eqref{eq:rbl3}''
and \eqref{adrvb3a}, respectively.}
\label{fig1b}
\end{figure}

In Figure \ref{fig1b}, we implement $6$ outer steps
of the restarted block Lanczos method \eqref{rbl}
with the block size $3$. Each outer step contains $4$ inner steps.
The comparison again utilizes mean values
among $1000$ samples with random initial subspaces.
Mean values of Ritz value errors from \eqref{rbl} and its stepwise modification
by the shifted Chebyshev polynomial \eqref{fch} are contained in
``Lanczos'' curves and ``Chebyshev'' circles. We count
$4{\,+\,}(4{\,-\,}1){\,\times\,}5{\,=\,}19$ inner steps. Restart occurs
in the iteration indices $4,\,7,\,10,\,13,\,16$ where we use
Theorem \ref{thm:rblm} for determining several nodes
in ``$\mathrm{Bound}_1$'' and ``$\mathrm{Bound}_2$''.
The other nodes concerning the rest inner steps are generated
by Theorem \ref{thm:rbl} as above for Figure \ref{fig1a}.
Thus $\mathrm{Bound}_1$ illustrates a combination
``\eqref{eq:rblm2}$+$\eqref{eq:rbl2}'', and $\mathrm{Bound}_2$
corresponds to ``\eqref{eq:rblm3}$+$\eqref{eq:rbl3}''.
In addition, based on \eqref{adrvb3}, we evaluate
\begin{equation}\label{adrvb3a}
 \left[T_{c}\left(1+2\,\frac{\la_i-\la_{p+1}}{\la_{p+1}-\la_n}\right)\right]^{-2}
 \left[T_{k}\left(1+2\,\frac{\la_i-\la_{p+1}}{\la_{p+1}-\la_n}\right)\right]^{-2(s-1)}
 \tan^2\angle(\X_i,\Y_i^{(\ell)})
\end{equation}
in the $(c{\,+\,}1)$th inner step of the $s$th outer step for generating ``$\mathrm{Bound}_3$''.
These three bounds behave similar to their counterparts in Figure \ref{fig1a}.
$\mathrm{Bound}_1$ generally improves $\mathrm{Bound}_2$
and $\mathrm{Bound}_3$ by using nonconsecutive eigenvalues
in convergence factors and angle-free constant terms.

\subsection{Example 2}

We reuse the test matrix from \cite[Example 2]{zakn2022}, i.e.,
$A=\mbox{diag}(\la_1,\ldots,\la_n)$ with $n=3600$ and
\begin{equation*}\begin{split}
 &\la_1=2.05, \quad \la_2=2, \quad \la_3=1.95, \quad
   \la_4=1.65, \quad \la_5=1.6, \quad \la_6=1.55,\\
 &\la_7=1.45, \quad \la_8=1.4, \quad \la_9=1.35, \quad
   \la_j=1-(j-9)/n \ \ \mbox{for} \ \ j=10,\ldots,n.
\end{split}\end{equation*}
Therein $\la_1,\ldots,\la_9$ are considered as target eigenvalues
and build three clusters. Correspondingly, we implement
the restarted block Lanczos method \eqref{rbl} for the block size $p=9$.
Further settings and bound evaluations are the same as those used in Example 1.

\begin{figure}[htbp]
\begin{center}
\includegraphics[width=\textwidth]{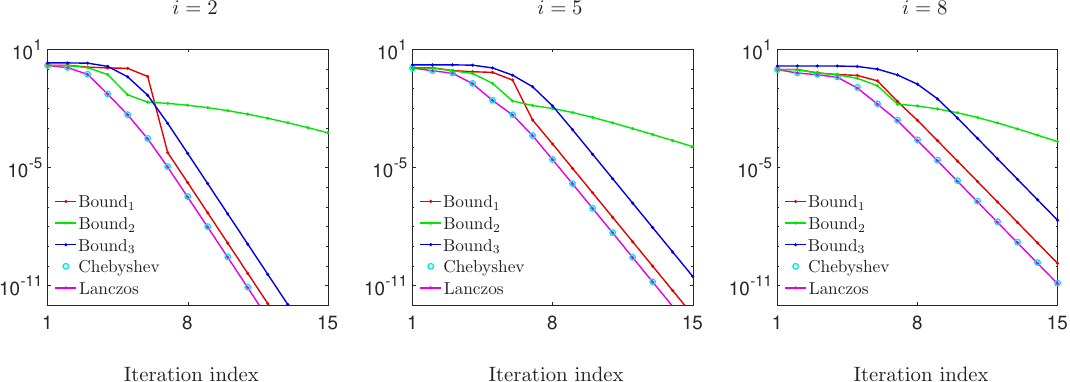}\quad
\end{center}
\par\vskip -2ex
\caption{\small Numerical comparison between several Ritz value bounds
concerning one outer step of the restarted block Lanczos method \eqref{rbl}
in Example 2. The three bound curves are determined by bounds
\eqref{eq:rbl2}, \eqref{eq:rbl3} and \eqref{adrvb3}, respectively.}
\label{fig2a}
\end{figure}

\begin{figure}[htbp]
\begin{center}
\includegraphics[width=\textwidth]{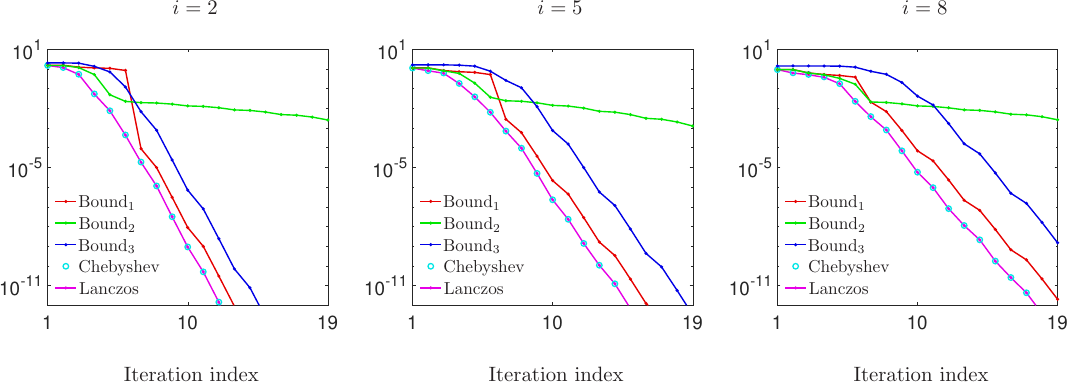}\quad
\end{center}
\par\vskip -2ex
\caption{\small Numerical comparison between several Ritz value bounds
concerning $6$ outer steps ($19$ inner steps)
of the restarted block Lanczos method \eqref{rbl}
in Example 2. The three bound curves are determined by bound combinations
``\eqref{eq:rblm2}$+$\eqref{eq:rbl2}'', ``\eqref{eq:rblm3}$+$\eqref{eq:rbl3}''
and \eqref{adrvb3a}, respectively.}
\label{fig2b}
\end{figure}

Figure \ref{fig2a} and Figure \ref{fig2b} illustrate one outer step ($15$ inner steps)
and $6$ outer steps ($19$ inner steps) of \eqref{rbl}, respectively. We document
Ritz value errors for the indices $2,\,5,\,8$ regarding the eigenvalue clusters.
The cluster robustness of \eqref{rbl} is clearly reflected by
$\mathrm{Bound}_1$ and $\mathrm{Bound}_3$, whereas
$\mathrm{Bound}_2$ depending on consecutive eigenvalues
gives a more substantial overestimation in the final phase
in comparison to Example 1. Moreover, the distances between
$\mathrm{Bound}_1$ and $\mathrm{Bound}_3$ become more evident
due to relatively larger angle terms.

\section*{Conclusions}

The convergence theory of block eigensolvers \cite{s1980,ks,lz2015}
is dominated by angle-dependent bounds where angles
build constant terms beside central convergence factors.
For analyzing restarted versions, it is necessary to adjust these terms
for associated convergence measures. The consequent overestimation is expected
to be avoided by direct bounds which can be applied repeatedly.
We extend one such bound, the angle-free Ritz value bound \cite[(2.22)]{ks},
to more Ritz values and flexible locations. Therein an argument
with biorthogonal vectors for investigating the block power method
by Rutishauser \cite{r1969} is modified by considering intersections
of iterative subspaces and certain invariant subspaces. This enables
generating partial iterations with which the derivation of \cite[(2.22)]{ks}
is upgraded. The resulting bounds are particularly applicable
to the restarted block Lanczos method for multiple outer steps,
and improve comparable angle-dependent bounds.
An open issue is whether the above approach can be generalized to
sums of Ritz value errors for supplementing angle-dependent
majorization bounds from \cite{zakn2022}. Furthermore, the application
to restarted block eigensolvers with shift-and-invert is potentially extendable to
corresponding preconditioned versions, similarly to recent results
on block preconditioned gradient-type eigensolvers in \cite{zn2019,zn2022}.

\section{Appendix}

\subsection{Proof of Lemma \ref{lm:afrvb1}}\label{pr1}

The assumption $\eta_p>\la_{p+1}$ ensures that
$X^HY$ has full rank for $X=[x_1,\ldots,x_p]$,
since otherwise there exists a vector $\widehat{g}\in\C^p{\setminus}\{0\}$
satisfying $X^HY\widehat{g}=0$ and causing a contradiction:
The vector $Y\widehat{g}\in\Y$ belongs to the invariant subspace
$\span\{x_{p+1},\ldots,x_n\}$ so that $\la_{p+1}\ge\rho(Y\widehat{g})\ge\eta_p$
holds and contradicts $\eta_p>\la_{p+1}$.

Moreover, the diagonal matrix $D=\diag\big(f(\la_1),\ldots,f(\la_p)\big)$
is invertible due to the assumption on $f(\cdot)$ so that $DX^HY$
has full rank, and $f(A)Y$ also does since $DX^HY=X^Hf(A)Y$. Then
$f(A)Y$ is a basis matrix of $\Y'$ since $\Y'=f(A)\,\span\{Y\}=\span\{f(A)Y\}$.
This trivially implies $\dim\Y'=p$.
\hfill$\Box$

\subsection{Proof of Lemma \ref{lm:afrvb2}}\label{app1}

The relations $\la_p>\eta'_p=\rho(y')$ and $\rho(y)\ge\eta_p>\la_{p+1}$
in \eqref{afrvbr} simply follow from the settings and the Courant-Fischer principles.
The remaining relation $\rho(y')\ge\rho(y^{\circ})\ge\rho(y)$ can be shown
in three steps by Lemma \ref{lm:rqw}.

(i) The relation $\rho(y)\ge\eta_p>\la_{p+1}$
ensures $\la_l\ge\rho(y)\ge\la_{l+1}$ for a certain index $l\in\{1,\ldots,p\}$.
According to $|f(\la_1)|\ge\cdots\ge|f(\la_p)|>\max_{j\in\{p+1,\ldots,n\}}|f(\la_j)|\ge0$,
it holds that
\[\frac{|f(\la_1)|}{|f(\la_l)|}\ge\cdots\ge\frac{|f(\la_l)|}{|f(\la_l)|}=1
 \ge\cdots\ge\frac{|f(\la_p)|}{|f(\la_l)|},\]
and
\[\frac{|f(\la_j)|}{|f(\la_l)|}\le\frac{\max_{j\in\{p+1,\ldots,n\}}|f(\la_j)|}
 {|f(\la_p)|}<1 \ \ \forall\ j>p.\]
By defining $w'=y'/f(\la_l)=f(A)y/f(\la_l)$, we get $\rho(w')=\rho(y')$, and
\[|x_j^Hw'|=\frac{|x_j^Hf(A)y|}{|f(\la_l)|}
 =\frac{|f(\la_j)x_j^Hy|}{|f(\la_l)|}
 =\frac{|f(\la_j)|}{|f(\la_l)|}\,|x_j^Hy|
 \left\{\begin{array}{ll}
 \ge |x_j^Hy| & \forall\ j \le l, \\[1ex]
 \le |x_j^Hy| & \forall\ j > l.
 \end{array}\right.\]
Then applying Lemma \ref{lm:rqw} (a) to
$u=y$ and $v=w'$ yields $\rho(w')\ge\rho(y)$ so that
\[\la_p>\eta'_p=\rho(y')=\rho(w')\ge\rho(y)>\la_{p+1}.\]

(ii) The assumption on $f(\cdot)$ ensures $f(\la_p)\neq0$.
Then $y^{\circ}\neq0$ holds since otherwise
\[0=x_j^Hy^{\circ}=f(\la_p)x_j^Hy \quad \forall \ \ j \le p
 \ \quad\Rightarrow\quad \
 x_j^Hy=0 \quad \forall \ \ j \le p
 \ \quad\Rightarrow\quad \
 \rho(y)\le\la_{p+1} \]
contradicts $\rho(y)\ge\eta_p>\la_{p+1}$.
Thus $\mu(y^{\circ})$ can be defined.
According to the result $\la_p>\rho(y')>\la_{p+1}$ from (i) and
\[|x_j^Hy^{\circ}|=\left\{\begin{array}{ll}
 |f(\la_p)|\,|x_j^Hy|\le|f(\la_j)|\,|x_j^Hy|
 =|x_j^Hf(A)y|=|x_j^Hy'| & \forall\ j \le p, \\[1ex]
 \big(\max_{j\in\{p+1,\ldots,n\}}|f(\la_j)|\big)\,|x_j^Hy|
 \ge|f(\la_j)|\,|x_j^Hy|=|x_j^Hy'| & \forall\ j > p,
 \end{array}\right.\]
we can apply Lemma \ref{lm:rqw} (b) to $l=p$,
$u=y'$ and $v=y^{\circ}$. This implies $\rho(y')\ge\rho(y^{\circ})$.

(iii) The vector $w^{\circ}=y^{\circ}/f(\la_p)\neq0$ fulfills
\[|x_j^Hw^{\circ}|=\left\{\begin{array}{ll}
 |f(\la_p)|\,|x_j^Hy| / |f(\la_p)|=|x_j^Hy| & \forall\ j \le p, \\[1ex]
 \big(\max_{j\in\{p+1,\ldots,n\}}|f(\la_j)|\big)\,|x_j^Hy| / |f(\la_p)|
 \le|x_j^Hy| & \forall\ j > p.
 \end{array}\right.\]
Using this together with the result $\la_p>\rho(y)>\la_{p+1}$ from (i) 
allows applying Lemma \ref{lm:rqw} (a) to $l=p$,
$u=y$ and $v=w^{\circ}$. Then we get
$\rho(y^{\circ})=\rho(w^{\circ})\ge\rho(y)$.

Summarizing (ii) and (iii) yields $\rho(y')\ge\rho(y^{\circ})\ge\rho(y)$.

For deriving \eqref{afrvbs}, we denote $y_p$ and $y-y_p$
by $w$ and $z$ so that
\[y=w+z \quad\mbox{and}\quad y^{\circ}=f(\la_p)w+\nu_p z.\]
Then $y$ and $y^{\circ}$ can be analyzed within the subspace $\span\{w,\,z\}$.

Therein $w$ and $z$ are nonzero vectors since otherwise
\,$y=z$\, or \,$y=w$\, holds so that \,$\rho(y)\le\la_{p+1}$\,
or \,$\rho(y)\ge\la_p$\, contradicts \eqref{afrvbr}.
Thus $\rho(w)$ and $\rho(z)$ can be defined. Applying the orthogonality
properties \,$w^Hz=0$\, and \,$w^HAz=0$\, to \,$y=w+z$\, gives
\begin{equation}\label{eq:app1}
\begin{split}
 &\rho(y)=\frac{y^HAy}{y^Hy}=\frac{w^HAw+z^HAz}{w^Hw+z^Hz}
 =\frac{\rho(w)\|w\|_2^2+\rho(z)\|z\|_2^2}{\|w\|_2^2+\|z\|_2^2}
 \\&\hspace{2cm}
 \quad\Rightarrow\quad
 \frac{\rho(w)-\rho(y)}{\rho(y)-\rho(z)}=\frac{\|z\|_2^2}{\|w\|_2^2}.
\end{split}
\end{equation}
Combining this with an analogous result
for \,$y^{\circ}=f(\la_p)w+\nu_p z$\, yields
\[\frac{\rho(w)-\rho(y^{\circ})}{\rho(y^{\circ})-\rho(z)}
 =\frac{\|\nu_p z\|_2^2}{\|f(\la_p)w\|_2^2}
 =\frac{\nu_p^2}{|f(\la_p)|^2}\,\frac{\|z\|_2^2}{\|w\|_2^2}
 =\frac{\nu_p^2}{|f(\la_p)|^2}\,\frac{\rho(w)-\rho(y)}{\rho(y)-\rho(z)}\]
and implies \eqref{afrvbs}.
\hfill$\Box$

\subsection{Proof of Theorem \ref{thm:afrvb}}\label{pr2}

The assumption on $f(\cdot)$ ensures $f(\la_j)\neq0$ for $j \le p$
so that Lemma \ref{lm:afrvb1} is applicable and gives $\dim\Y'=p$.
For verifying \eqref{afrvb}, we can skip the trivial case $\la_p=\eta'_p$.
Then $\la_p>\eta'_p$ holds so that Lemma \ref{lm:afrvb2} is applicable.
Denoting $y_p$ and $y-y_p$ by $w$ and $z$, \eqref{afrvbr} implies
\[\rho(w)\ge\la_p>\eta'_p=\rho(y')\ge\rho(y^{\circ})\ge\rho(y)
 \ge\eta_p>\la_{p+1}\ge\rho(z).\]
Subsequently, simple monotonicity arguments lead to
\begin{equation}\label{eq:afrvb}
\begin{split}
 &\left(\frac{\la_p-\eta'_p}{\eta'_p-\la_{p+1}}\right)
 \left(\frac{\la_p-\eta_p}{\eta_p-\la_{p+1}}\right)^{-1}\le
 \left(\frac{\la_p-\rho(y^{\circ})}{\rho(y^{\circ})-\la_{p+1}}\right)
 \left(\frac{\la_p-\rho(y)}{\rho(y)-\la_{p+1}}\right)^{-1} \\[1ex]
 =\ &\left(\frac{\la_p-\rho(y^{\circ})}{\la_p-\rho(y)}\right)
 \left(\frac{\rho(y)-\la_{p+1}}{\rho(y^{\circ})-\la_{p+1}}\right)
 \le\left(\frac{\rho(w)-\rho(y^{\circ})}{\rho(w)-\rho(y)}\right)
 \left(\frac{\rho(y)-\rho(z)}{\rho(y^{\circ})-\rho(z)}\right) \\[1ex]
 =\ &\left(\frac{\rho(w)-\rho(y^{\circ})}{\rho(y^{\circ})-\rho(z)}\right)
 \left(\frac{\rho(w)-\rho(y)}{\rho(y)-\rho(z)}\right)^{-1}
 \stackrel{\eqref{afrvbs}}{=}\frac{\max_{j\in\{p+1,\ldots,n\}}|f(\la_j)|^2}{|f(\la_p)|^2}
\end{split}
\end{equation}
which results in \eqref{afrvb}.
\hfill$\Box$

\subsection{Proof of Lemma \ref{lm:aux}}\label{pr3}

We denote the three subspaces with the subscript ${[1,i]\cup(p,n]}$
by $\widetilde{\X}$, $\widetilde{\Y}$ and $\widetilde{\Y}'$, respectively.
The dimension comparison
\[\dim(\widetilde{\X}\cap\Y)=\dim\widetilde{\X}
 +\dim\Y-\dim(\widetilde{\X}+\Y)\ge (n-p+i)+p-n = i\]
gives $\dim\widetilde{\Y} \ge i$.
In the case $\eta_p>\la_{p+1}$, the strict inequality $\dim\widetilde{\Y}>i$
does not hold, since otherwise the smallest 
($j$th largest with $j>i$) Ritz value $\widetilde{\eta}$
in $\widetilde{\Y}$ fulfills $\widetilde{\eta}\le\la_{p+1}$
and $\widetilde{\eta}\ge\eta_p$ 
due to $\widetilde{\Y}\subseteq\widetilde{\X}$
and $\widetilde{\Y}\subseteq\Y$ so that $\eta_p\le\la_{p+1}$.

As shown in the proof of Lemma \ref{lm:afrvb1},
$X^HY$ is an invertible $p{\,\times\,}p$ matrix for $X=[x_1,\ldots,x_p]$
and an arbitrary basis matrix $Y$ of $\Y$. 
Then the $n{\,\times\,}i$ matrix $Y_i=Y(X^HY)^{-1}[e_1,\ldots,e_i]$ has full rank
where $e_1,\ldots,e_p$ are columns of the $p{\,\times\,}p$ identity matrix.
Thus $\span\{Y_i\}$ is an $i$-dimensional subset of $\Y$. Moreover, since
$[x_{i+1},\ldots,x_p]^HY_i$ is a zero matrix due to
\begin{equation*}\begin{split}
 [x_{i+1},\ldots,x_p]^HY_i
 &=[e_{i+1},\ldots,e_p]^HX^HY(X^HY)^{-1}[e_1,\ldots,e_i]\\
 &=[e_{i+1},\ldots,e_p]^H[e_1,\ldots,e_i],
\end{split}\end{equation*}
it holds that $\span\{Y_i\}\subseteq\widetilde{\X}$ and consequently
$\span\{Y_i\}=\widetilde{\Y}$. By using $XX^H$ as the orthogonal projector
$P_{\X}$ on $\X$, we get
\[(XX^H)Y_i=XX^HY(X^HY)^{-1}[e_1,\ldots,e_i]=X\,[e_1,\ldots,e_i]=[x_1,\ldots,x_i]\]
so that $P_{\X}\widetilde{\Y}=\span\{x_1,\ldots,x_i\}$.

In addition, for $X_i=[x_1,\ldots,x_i]$ and $D_i=\diag\big(f(\la_1),\ldots,f(\la_i)\big)$,
it holds that $X_i^Hf(A)Y_i=D_iX_i^HY_i=D_i$, and 
the assumption on $f(\cdot)$ ensures that $D_i$ is invertible. 
Thus $f(A)Y_i$ has full rank so that the subspace
$\widetilde{\Y}'=f(A)\,\span\{Y_i\}=\span\{f(A)Y_i\}$ has dimension $i$.
\hfill$\Box$

\subsection{Proof of Lemma \ref{lm:afrvb3}}\label{app2}

Following the proof of Lemma \ref{lm:aux} with simplified subspace notation,
we use again the basis matrices $Y_i$ and $f(A)Y_i$
of $\widetilde{\Y}$ and $\widetilde{\Y}'$. Then an arbitrary Ritz vector $\widetilde{y}'$
associated with $\widetilde{\eta}'_i$ can be represented by
$\widetilde{y}'=f(A)Y_ig$ with a certain $g\in\C^i{\setminus}\{0\}$
so that $\widetilde{y}'=f(A)\widetilde{y}$ for $\widetilde{y}=Y_ig$.

For verifying \eqref{afrvbr1},
the relation $\la_i>\widetilde{\eta}'_i=\rho(\widetilde{y}')$
is trivial. Moreover, $\widetilde{y}\in\widetilde{\Y}\subseteq\Y$ implies
$\rho(\widetilde{y})\ge\widetilde{\eta}_i\ge\eta_p$. Combining this
with the assumption $\eta_p>\la_{p+1}$ from Lemma \ref{lm:aux}
gives $\rho(\widetilde{y})\ge\widetilde{\eta}_i>\la_{p+1}$.
Subsequently, $\rho(\widetilde{y}')\ge\rho(\widetilde{y}^{\circ})
\ge\rho(\widetilde{y})$ follows from Lemma \ref{lm:rqw}
in three steps, with some nontrivial detailed changes
in comparison to the proof of Lemma \ref{lm:afrvb2}.

(i) We have $\la_l\ge\rho(\widetilde{y})\ge\la_{l+1}$
for a certain index $l\in\{1,\ldots,p\}$ due to
$\rho(\widetilde{y})\ge\widetilde{\eta}_i>\la_{p+1}$.
If $l < i$, the assumption
$|f(\la_1)|\ge\cdots\ge|f(\la_i)|>\max_{j\in\{p+1,\ldots,n\}}|f(\la_j)|$ leads to
\[\frac{|f(\la_1)|}{|f(\la_l)|}\ge\cdots\ge\frac{|f(\la_l)|}{|f(\la_l)|}=1
 \ge\cdots\ge\frac{|f(\la_i)|}{|f(\la_l)|},\]
and
\[\frac{|f(\la_j)|}{|f(\la_l)|}\le\frac{\max_{j\in\{p+1,\ldots,n\}}|f(\la_j)|}
 {|f(\la_i)|}<1 \ \ \forall\ j>p.\]
Then $\widetilde{w}'=\widetilde{y}'/f(\la_l)=f(A)\widetilde{y}/f(\la_l)$
fulfills $\rho(\widetilde{w}')=\rho(\widetilde{y}')$,
\[|x_j^H\widetilde{w}'|=\frac{|x_j^Hf(A)\widetilde{y}|}{|f(\la_l)|}
 =\frac{|f(\la_j)x_j^H\widetilde{y}|}{|f(\la_l)|}
 \left\{\begin{array}{ll}
 \ge |x_j^H\widetilde{y}| & \forall\ j \le l, \\[1ex]
 \le |x_j^H\widetilde{y}| & \forall\ j \in\{l{\,+\,}1,\ldots,i,\,p{\,+\,}1,\ldots,n\}.
 \end{array}\right.\]
If $l \ge i$, we redefine $\widetilde{w}'$ by
$\widetilde{w}'=\widetilde{y}'/f(\la_i)$ so that
\[|x_j^H\widetilde{w}'|=\frac{|x_j^Hf(A)\widetilde{y}|}{|f(\la_i)|}
 =\frac{|f(\la_j)x_j^H\widetilde{y}|}{|f(\la_i)|}
 \left\{\begin{array}{ll}
 \ge |x_j^H\widetilde{y}| & \forall\ j \le i, \\[1ex]
 \le |x_j^H\widetilde{y}| & \forall\ j >p.
 \end{array}\right.\]
In both cases, a combination with
$|x_j^H\widetilde{w}'|=0=|x_j^H\widetilde{y}|
\ \ \forall \ j \in\{i{\,+\,}1,\ldots,p\}$ gives
\[|x_j^H\widetilde{w}'|\ge|x_j^H\widetilde{y}|
\ \ \forall \ j \le l
\quad\mbox{and}\quad
|x_j^H\widetilde{w}'|\le|x_j^H\widetilde{y}|
\ \ \forall \ j > l.\]
Thus Lemma \ref{lm:rqw} (a) is applicable to
$u=\widetilde{y}$ and $v=\widetilde{w}'$, and implies
\[\la_i>\eta'_i=\rho(\widetilde{y}')=\rho(\widetilde{w}')
 \ge\rho(\widetilde{y})>\la_{p+1}.\] 

(ii) The vector $\widetilde{y}^{\circ}$ is nonzero since otherwise
\[0=x_j^H\widetilde{y}^{\circ}=f(\la_i)x_j^H\widetilde{y} \quad \forall \ \ j \le i
 \quad\Rightarrow\quad
 x_j^H\widetilde{y}=0 \quad \forall \ \ j \le i
 \quad\Rightarrow\quad
 x_j^H\widetilde{y}=0 \quad \forall \ \ j \le p\]
so that $\rho(\widetilde{y})\le\la_{p+1}$ holds
and contradicts $\rho(\widetilde{y})\ge\widetilde{\eta}_i>\la_{p+1}$.
Thus $\mu(y^{\circ})$ can be defined.
The result $\la_i>\rho(\widetilde{y}')>\la_{p+1}$ from (i)
indicates $\la_l\ge\rho(\widetilde{y}')\ge\la_{l+1}$
for a certain index $l\in\{i,\ldots,p\}$. In addition, combining
\[|x_j^H\widetilde{y}^{\circ}|=\left\{\begin{array}{ll}
 |f(\la_i)|\,|x_j^H\widetilde{y}|\le|f(\la_j)|\,|x_j^H\widetilde{y}|
 =|x_j^Hf(A)\widetilde{y}|=|x_j^H\widetilde{y}'| & \forall\ j \le i, \\[1ex]
 \big(\max_{j\in\{p+1,\ldots,n\}}|f(\la_j)|\big)\,|x_j^H\widetilde{y}|
 \ge|f(\la_j)|\,|x_j^H\widetilde{y}|=|x_j^H\widetilde{y}'| & \forall\ j > p
 \end{array}\right.\]
with $|x_j^H\widetilde{y}^{\circ}|=0=|x_j^H\widetilde{y}'|
\ \ \forall \ j \in\{i{\,+\,}1,\ldots,p\}$ ensures that Lemma \ref{lm:rqw} (b)
is applicable to $u=\widetilde{y}'$ and $v=\widetilde{y}^{\circ}$.
Thus $\rho(\widetilde{y}')\ge\rho(\widetilde{y}^{\circ})$.

(iii) It holds that $\la_l\ge\rho(\widetilde{y})\ge\la_{l+1}$
for a certain index $l\in\{i,\ldots,p\}$
due to $\la_i>\rho(\widetilde{y})>\la_{p+1}$ from (i).
The auxiliary vector $\widetilde{w}^{\circ}=\widetilde{y}^{\circ}/f(\la_i)$
is nonzero and fulfills
\[|x_j^H\widetilde{w}^{\circ}|=\left\{\begin{array}{ll}
 |f(\la_i)|\,|x_j^H\widetilde{y}| / |f(\la_i)|
 =|x_j^H\widetilde{y}| & \forall\ j \le i, \\[1ex]
 \big(\max_{j\in\{p+1,\ldots,n\}}|f(\la_j)|\big)\,|x_j^H\widetilde{y}| / |f(\la_i)|
 \le|x_j^H\widetilde{y}| & \forall\ j > p
 \end{array}\right.\]
and $|x_j^H\widetilde{w}^{\circ}|=0=|x_j^H\widetilde{y}|
\ \ \forall \ j \in\{i{\,+\,}1,\ldots,p\}$. This allows applying Lemma \ref{lm:rqw} (a)
to $u=\widetilde{y}$ and $v=\widetilde{w}^{\circ}$ so that
$\rho(\widetilde{y}^{\circ})=\rho(\widetilde{w}^{\circ})\ge\rho(\widetilde{y})$.

According to (ii) and (iii), we get
$\rho(\widetilde{y}')\ge\rho(\widetilde{y}^{\circ})\ge\rho(\widetilde{y})$.

The derivation of \eqref{afrvbs1} is based on the representations
\[\widetilde{y}=\widetilde{w}+\widetilde{z} \quad\mbox{and}\quad
 \widetilde{y}^{\circ}=f(\la_i)\widetilde{w}+\nu_p \widetilde{z}\]
with $\widetilde{w}=\widetilde{y}_i$ and $\widetilde{z}=\widetilde{y}-\widetilde{y}_i$.
The property \eqref{afrvbr1} excludes $\widetilde{w}=0$
or $\widetilde{z}=0$ which would lead to 
\,$\rho(\widetilde{y})\le\la_{p+1}$\, or \,$\rho(\widetilde{y})\ge\la_i$.
Thus $\rho(\widetilde{w})$ and $\rho(\widetilde{z})$ can be defined.
By using \,$\widetilde{w}^H\widetilde{z}=0$\,
and \,$\widetilde{w}^HA\widetilde{z}=0$,
\[\frac{\rho(\widetilde{w})-\rho(\widetilde{y})}
 {\rho(\widetilde{y})-\rho(\widetilde{z})}
 =\frac{\|\widetilde{z}\|_2^2}{\|\widetilde{w}\|_2^2}
 \quad\mbox{and}\quad
 \frac{\rho(\widetilde{w})-\rho(\widetilde{y}^{\circ})}
 {\rho(\widetilde{y}^{\circ})-\rho(\widetilde{z})}
 =\frac{\|\nu_p \widetilde{z}\|_2^2}{\|f(\la_i)\widetilde{w}\|_2^2}\]
can be shown analogously to \eqref{eq:app1}, and result in \eqref{afrvbs1}.
\hfill$\Box$

\subsection{Proof of Lemma \ref{lm:afrvb4}}\label{app3}

The property $\dim\widetilde{\Y} \ge i$ follows from a dimension comparison
as in the proof of Lemma \ref{lm:aux}.

For verifying $\dim\widetilde{\Y}'_i=i$, we use an arbitrary basis matrix
$\widetilde{Y}_i\in\C^{n \times i}$ of $\widetilde{\Y}_i$
together with $X_{t-p+i}=[x_1,\ldots,x_{t-p+i}]$. The assumption
$\eta_p>\la_{t+1}$ ensures that $X_{t-p+i}^T\widetilde{Y}_i$
has full rank since otherwise there exists a vector $\widehat{g}\in\C^i{\setminus}\{0\}$
with $X_{t-p+i}^T\widetilde{Y}_i\widehat{g}=0$ so that
the vector $\widehat{y}=\widetilde{Y}_i\widehat{g}$
is orthogonal to $x_1,\ldots,x_{t-p+i}$. Combining this with
$\widehat{y}\in\widetilde{\Y}_i\subseteq\widetilde{\Y}
\subseteq\widetilde{\X}=\span\{x_{t-p+i+1},\ldots,x_t\}^{\perp}$
shows that $\widehat{y}$ belongs to $\span\{x_{t+1},\ldots,x_n\}$.
Then $\la_{t+1}\ge\rho(\widehat{y})\ge\eta_p$ holds
(due to $\widehat{y}\in\widetilde{\Y}_i\subseteq\widetilde{\Y}\subseteq\Y$)
and contradicts $\eta_p>\la_{t+1}$. In addition, the assumption
$f(\la_j)\neq0 \ \ \forall \ j\in\{1,\ldots,t{\,-\,}p{\,+\,}i\}$
suggests the invertible diagonal matrix
$D_{t-p+i}=\diag\big(f(\la_1),\ldots,f(\la_{t-p+i})\big)$.
Therewith $X_{t-p+i}^Tf(A)\widetilde{Y}_i=D_{t-p+i}X_{t-p+i}^T\widetilde{Y}_i$
has full rank so that $f(A)\widetilde{Y}_i$ also does,
and $\widetilde{\Y}'_i=f(A)\,\span\{\widetilde{Y}_i\}
=\span\{f(A)\widetilde{Y}_i\}$ has dimension $i$.

Furthermore, by using $f(A)\widetilde{Y}_i$ as a basis matrix of $\widetilde{\Y}'_i$,
we get the representation $\widetilde{y}'=f(A)\widetilde{Y}_ig$
with a certain $g\in\C^i{\setminus}\{0\}$. Thus
$\widetilde{y}'=f(A)\widetilde{y}$ for $\widetilde{y}=\widetilde{Y}_ig$.

The properties \eqref{afrvbr2} and \eqref{afrvbs2} can be shown by modifying
the derivations of \eqref{afrvbr1} and \eqref{afrvbs1} in Subsection \ref{app2}
with shifted indices.
\hfill$\Box$

\subsection{Strictly increasing Ritz values}\label{app4}

We explain the fact that the Ritz values strictly increase
during the restarted block Lanczos method \eqref{rbl}
as long as the current iterative subspace contains no eigenvectors.

With the settings from Theorem \ref{thm:rbl},
we show $\psi_i^{(\ell+1)}>\psi_i^{(\ell)}$ for each $i\in\{1,\ldots,p\}$
provided that there are no eigenvectors in $\Y^{(\ell)}$.
Therein we denote by $y_1,\ldots,y_p$ orthonormal Ritz vectors associated with
the Ritz values $\psi_1^{(\ell)}\ge\cdots\ge\psi_p^{(\ell)}$ in $\Y^{(\ell)}$,
and define $\widetilde{\Y}_i=\span\{y_1,\ldots,y_i\}$.
By using an arbitrary $\la<\la_n$, the matrix $\widetilde{A}=A-\la I$
is Hermitian positive definite, and the subspace
$\widetilde{\Y}'_i=\widetilde{A}\widetilde{\Y}_i$ has dimension $i$.
Then the $i$th largest Ritz value $\widetilde{\eta}'_i$ in $\widetilde{\Y}'_i$
fulfills $\widetilde{\eta}'_i\le\psi_i^{(\ell+1)}$ according to
\[\widetilde{A}\widetilde{\Y}_i\subseteq\widetilde{A}\Y^{(\ell)}
 \subseteq\Y^{(\ell)}+A\Y^{(\ell)}+\cdots+A^{k}\Y^{(\ell)}.\]
Moreover, an arbitrary Ritz vector $\widetilde{y}'$
associated with $\widetilde{\eta}'_i$ can be represented by
\,$\widetilde{y}'=\widetilde{A}\widetilde{y}$\,
with a certain $\widetilde{y}\in\widetilde{\Y}_i{\setminus}\{0\}$.

We further use the Rayleigh quotients $\rho(\cdot)$ and $\widetilde{\rho}(\cdot)$
with respect to $A$ and $\widetilde{A}$. Then
\[\widetilde{\rho}(\widetilde{y})
 =\frac{\widetilde{y}^H\widetilde{A}\widetilde{y}}{\widetilde{y}^H\widetilde{y}}
 \le\frac{\|\widetilde{y}\|_2\|\widetilde{A}\widetilde{y}\|_2}
 {\|\widetilde{y}\|_2^2}
 =\frac{(\widetilde{y},\,\widetilde{A}\widetilde{y})_{\widetilde{A}}}
 {\|\widetilde{y}\|_2\|\widetilde{A}\widetilde{y}\|_2}
 \le\frac{\|\widetilde{y}\|_{\widetilde{A}}\|\widetilde{A}\widetilde{y}\|_{\widetilde{A}}}
 {\|\widetilde{y}\|_2\|\widetilde{A}\widetilde{y}\|_2}
 =\big(\widetilde{\rho}(\widetilde{y})\,
 \widetilde{\rho}(\widetilde{y}')\big)^{1/2}\]
implies $\widetilde{\rho}(\widetilde{y})\le\widetilde{\rho}(\widetilde{y}')$,
i.e., $\rho(\widetilde{y})-\la\le\rho(\widetilde{y}')-\la$
so that $\rho(\widetilde{y})\le\rho(\widetilde{y}')$.
Therein the equality holds if and only if $\widetilde{y}$ is collinear with
$\widetilde{A}\widetilde{y}$, i.e., $\widetilde{y}$ is an eigenvector of
$\widetilde{A}$ as well as $A$. Consequently, if $\Y^{(\ell)}$ contains
no eigenvectors, we get $\rho(\widetilde{y})<\rho(\widetilde{y}')$ so that
$\psi_i^{(\ell+1)}\ge\widetilde{\eta}'_i=\rho(\widetilde{y}')
>\rho(\widetilde{y})\ge\psi_i^{(\ell)}$.

\vskip12pt

\def\refname{
\end{document}